\documentclass[12pt]{article}
\usepackage{amsfonts}
\usepackage{amsmath}
\usepackage{pdfsync}

\usepackage[english]{babel}
\usepackage[latin1]{inputenc}
\usepackage[T1]{fontenc}
\usepackage{amsthm}
\usepackage[dvipsnames]{xcolor}

\newtheorem{theorem}{Theorem}[section]

\newtheorem{lemma}[theorem]{Lemma}
\newtheorem{remark}[theorem]{Remark}
\newtheorem{corollary}[theorem]{Corollary}
\newtheorem{proposition}[theorem]{Proposition}

\numberwithin{equation}{section}
\usepackage{lmodern}
\usepackage[colorlinks=true, allcolors=blue]{hyperref}
\usepackage[a4paper,top=3cm,bottom=2cm,left=3cm,right=3cm,marginparwidth=1.75cm]{geometry}

\begin{document}

\centerline{\large \textbf{The quantum mechanics canonically associated to}}
\centerline{\large \textbf{free probability}}
\centerline{\large \textbf{Part I: Free momentum and associated kinetic energy}}

\bigskip\bigskip

\centerline{\textbf{Luigi Accardi} }
\centerline{Centro Vito Volterra, Universit\`a di Roma "Tor Vergata"}
\centerline{Via Columbia, 2, 00133  Roma, Italy -- accardi@volterra.uniroma2.it}

\centerline{\textbf{Tarek Hamdi } }
\centerline{ Department of Management Information Systems, College of Business Management}
\centerline{Qassim University, Ar Rass, Saudi Arabia }
\centerline{and Laboratoire d'Analyse Math\'ematiques  et applications LR11ES11}
\centerline{Universit\'e de Tunis El Manar, Tunisie -- t.hamdi@qu.edu.sa}
\centerline{\textbf{Yun Gang Lu} }
\centerline{ Dipartimento di Matematica, Universit\`a di Bari}
\centerline{via Orabona, 4, 70125  Bari, Italy -- yungang.lu@uniba.it}

\bigskip\bigskip

\centerline{\textbf{Abstract}}

\noindent

\noindent
After a short review of the quantum mechanics canonically associated with a classical real valued
random variable with all moments, we begin to study the quantum mechanics canonically associated
to the \textbf{standard semi--circle random variable} $X$, characterized by the fact that its
probability distribution is the semi--circle law $\mu$ on $[-2,2]$.
We  prove that, in the identification of $L^2([-2,2],\mu)$ with the $1$--mode interacting Fock
space $\Gamma_{\mu}$, defined by the orthogonal polynomial gradation of $\mu$, $X$ is mapped
into position operator and its canonically associated momentum operator $P$ into $i$ times
the $\mu$--Hilbert transform $H_{\mu}$ on $L^2([-2,2],\mu)$.
In the first part of the present paper, after briefly describing the simpler case of the
$\mu$--harmonic oscillator, we find an explicit expression for the action,
on the $\mu$--orthogonal polynomials, of the semi--circle analogue of the translation group $e^{itP}$
and of the semi--circle analogue of the free evolution $e^{itP^2/2}$ respectively in terms of
Bessel functions of the first kind and of confluent hyper--geometric series.
These results require the solution of the \textit{inverse normal order problem} on the quantum algebra canonically associated to the classical semi--circle random variable and are derived
in the second part of the present paper.
Since the problem to determine, with purely analytic techniques, the explicit form of the action
of $e^{-tH_{\mu}}$ and $e^{-itH_{\mu}^2/2}$ on the $\mu$--orthogonal polynomials is difficult,
the above mentioned results show the power of the combination of these techniques with those
developed within the algebraic approach to the theory of orthogonal polynomials.

\tableofcontents

\section{Introduction}

\subsection{Short review of the quantum mechanics canonically associated with a classical random variable with all moments }

Let $X$ be a classical real valued random variable with all moments and probability
distribution $\mu$.
Denote $(\Phi_{n})_{n\in\mathbb{N}}$ the orthogonal polynomials of $X$, $(\omega_{n})$,
$(\alpha_{n})$ its Jacobi sequences and
\begin{equation}\label{1MIFS-om}
\Gamma_{X}\equiv\Gamma\left(  \mathbb{C},\left\{  \omega_{n}\right\}_{n=1}^{\infty}\right)
\subseteq L^2(\mathbb{R},\mu)
\end{equation}
the closure, in  $L^2(\mathbb{R},\mu)$, of the linear span of the $\Phi_{n}$'s.
It is known (see \cite{[AcBo98]}) that, identifying $X$ with the multiplication operator by $X$ in $\Gamma_{X}$,
the orthogonal gradation
\begin{equation}\label{1MIFS}
\Gamma_{X}\equiv \bigoplus_{n\in\mathbb{N}}\mathbb{C}\cdot\Phi_{n}
\end{equation}
uniquely defines, through Jacobi tri--diagonal relation, $3$ linear operators $ a^+, a^-, a^0$,
called respectively \textbf{creation, annihilation and preservation (CAP)} operators, leaving
invariant the linear span of the $\Phi_{n}$'s. In turn, the CAP operators uniquely define the
\textbf{canonical quantum decomposition} of $X$ through the identity
\begin{equation}\label{df-q-dec-X}
X =  a^+ + a^0 + a^-
\end{equation}
Moreover, denoting $\Lambda$ the unique linear extension, on the linear span of the $\Phi_{n}$'s,
of the map $\Phi_{n}\mapsto n\Phi_{n}$ (\textbf{number operator}) and defining, on the same
domain, for any sequence $(F_{n})$ of complex numbers, the linear operator $F_{\Lambda}$
by linear extension of $\Phi_{n}\mapsto F_{n}\Phi_{n}$, the operators  $ a^+, a^-$ and the
operator $\Lambda$ satisfy the following multiplication table
\begin{equation}\label{a+a-aa+}
a^+a^- = \omega_{\Lambda} \qquad;\qquad a^-a^+ = \omega_{\Lambda+1}
\end{equation}
\begin{equation}\label{F-Lambda-a}
F_{\Lambda}a^+=a^+F_{\Lambda+1}\qquad;\qquad a^-F_{\Lambda}=F_{\Lambda+1}a^-
\end{equation}
which implies the commutation relations
\begin{equation}\label{comm-rel-[a,a+]}
[a^-,a^+]=\omega_{\Lambda+1}-\omega_{\Lambda}=:\partial\omega_{\Lambda}
\qquad;\qquad[a^+,a^+]=[a^-,a^-]=0
\end{equation}
\begin{equation}\label{[a+,FLambda]}
[a^+,F_{\Lambda}]
= - a^+(F_{\Lambda+1}-F_{\Lambda})=: - a^+\partial F_{\Lambda}.
\end{equation}
The $*$--algebra generated by the CAP operators of $X$ and the functions of $\Lambda$ is called
the \textbf{quantum algebra canonically associated to the classical random variable $X$}.
The term \textbf{normal order} in this algebra is similar to that in usual quantum mechanics
with the only difference that, in this case,
the normally ordered expressions are sum of terms of the form
$(a^+)^ma^nF_{\Lambda}$, where $F_{\Lambda}$ is a function of $\Lambda$.\\
The Gauss and Poisson random variables have the same principal Jacobi sequence given by
$\omega_{n}=\hbar n$ where $\hbar $ is a strictly positive constant. Therefore, for
these measures, $\partial\omega_{\Lambda}=\hbar$ and the commutation relations
\eqref{comm-rel-[a,a+]} reduce to
\begin{equation}\label{Heis-CR}
[a^-,a^+]=\hbar
\qquad;\qquad[a^+,a^+]=[a^-,a^-]=0
\end{equation}
i.e. to the \textbf{Heisenberg commutation relations} for quantum systems with one degree of
freedom.\\
The decomposition \eqref{1MIFS} and the associated quantum decomposition of $X$ identifies
the theory of orthogonal polynomials in $1$ variable with the theory of
\textbf{$1$--mode interacting Fock space ($1$MIFS)} (see \cite{[AcLuVo97a-QED-Hilb-mod]}).\\
It is also known (see \cite{[AcEllLu20-mom]}) that the classical
random variable $X$ is (polynomially) \textbf{symmetric} (i.e. all its odd moments vanish)
if and only if in the decomposition \eqref{df-q-dec-X} $a^0=0$). For a symmetric classical
random variable $X$, the Hermitean operator
\begin{equation}\label{df-mom-op}
P_X:=i(a^+ - a)
\end{equation}
called the \textbf{momentum operator} canonically associated to $X$,
 satisfies the following commutation relation
\begin{equation}\label{[PX,X]}
[X,P_X] = i2\partial\omega_{\Lambda}
\end{equation}
which, for mean zero Gaussian random variables, reduces to the original Heisenberg
commutation relation between position and its canonical conjugate momentum.
Once the operators position $X$ and momentum $P_{X}$ are available, one can introduce
all operators of physical interest, like kinetic energy $P_{X}^2/2$, potential energy $\dots$.
Hence \textbf{any classical symmetric random variable with all moments uniquely determines its
own quantum mechanics} and usual quantum mechanics corresponds to mean zero Gaussian
random variables whose covariance plays the role of Planck's constant.\\

\noindent One knows that the centered Gaussian measure with variance $c$, is characterized by
$\omega_{n+1}-\omega_{n}=c$ for each $n\ge 0$, where $c$ is a \textbf{strictly positive} constant.
The case $\omega_{n+1}-\omega_{n}=0$ for each $n\ge 1$ (in section \ref{lm-semi-circle} it will
be clear why $1$ and not $0$) is mathematically and physically
intriguing because this case characterizes the principal Jacobi sequences of the form
\begin{equation}\label{semi-circ-princ-JS}
\omega_{n}=\omega \ge 0 \quad;\quad\forall n\in\mathbb{N}^*
\end{equation}
and it is known from classical probability that, if $\omega$ is a strictly positive constant,
the principal Jacobi sequences \eqref{semi-circ-princ-JS} \textbf{characterize the symmetric
semi--circle laws} with variance $\omega$. Since these laws play the role of the Gaussian for
free probability, the following problem naturally arises:
\textbf{which is the quantum mechanics canonically associated to free probability?}\\

\noindent In this paper $X$ will be the classical random variable (unique up to stochastic
equivalence) with probability distribution given by the semi--circle law with support in the
interval $(-2,2)$.
We prove (section \ref{analytic-forms}) that \textbf{the momentum operator $P_X$} canonically associated to $X$ is $i$ times \textbf{the
$\mu$--Hilbert transform} on $L^2([-2,2],\mu)$:
\begin{equation}\label{mu-Hilb_transf}
P_X = H_{\mu}f(x):= 2\hbox{p.v.}\int_{\mathbb{R}}\frac{f(y)}{x-y}d\mu(y)
\end{equation}
where $\hbox{p.v.}$ denotes Cauchy Principal Value.\\
Extensive research on the operator $H_{\mu}$ has been carried out both for its own interest
and for its important application to aerodynamics via the airfoil equation \cite{[King09]},
\cite{[Tricomi57]}.
Another interesting application of this operator includes tomography and problems arising in image reconstruction (see, e.g. \cite{[KaTo12]}).\\
This identification allows to combine the algebraic techniques developed in the IFS approach
to orthogonal polynomials with the analytical techniques developed in harmonic analysis to deal
with the Hilbert transform. As an illustration of this fact, we find the explicit form of the action,
on the monic orthogonal polynomials of $X$, of the $1$--parameter unitary groups generated by $P_X$
(section \ref{sec:evol-eitP}) and by the free kinetic energy operator $P_{X}^2/2$ associated to $X$, i.e. the free evolution (section \ref{sec:evol-eitP2}). \\
The solution of this problem requires the preliminary solution of the \textbf{inverse normal order
problem}. The normal and inverse normal order problems on the quantum algebra canonically
associated to the classical random variable $X$ are formulated and solved in the second part of
this paper.
These general results are applied in section
\ref{sec:The-ev-eitP-eitX} to deduce an explicit form for the Schr\"odinger and Heisenberg
evolutions associated to $e^{itP_{X}}$ and $e^{itX}$ and in section \ref{sec:evol-eitP2},  we solve
the same problem for $e^{itP_{X}^2}$ and $e^{itX^2}$.
In particular, we obtain an explicit study of the action of $ e^{-tH_\mu}$ and of $e^{-itH^2_\mu}$ on the $\mu$-orthogonal polynomials which is truly new and is a merit of our approach.

\subsection{The canonical quantum decomposition in the semi-circle case}\label{lm-semi-circle}

The unique symmetric classical real valued random variable with principal Jacobi sequence
$(\omega_{n})_{n\geq1}$ satisfying
\begin{equation}\label{om01}
\omega_{n}=\omega > 0
\quad,\quad \forall n\in\mathbb{N}^*:=\mathbb{N}\setminus\{0\}
\quad;\quad \omega_{0}:=0
\end{equation}
is called the \textbf{semi--circle random variable with parameter} $\omega$. Its law is
the semi--circle distribution with the same parameter and its canonical quantum decomposition is
$$
X:=a+a^{+}
$$
where $a, a^{+}$ are the creation--annihilation operators acting on the $1$MIFS
\eqref{1MIFS-om}, with $(\omega_{n})$ given by \eqref{om01}, and denoted
$$
\Gamma_{X}:=
\bigoplus_{n\in\mathbb{N}}\mathbb{C}\cdot\Phi_{n} = L^2(\mathbb{R},\mu).
$$
The monic basis of $\Gamma_{X}$ is
\begin{equation}\label{df-mon-bas}
\left\{\Phi_{n}:=a^{+n}\Phi_{0} \ : \  n\in\mathbb{N}\right\}
\end{equation}
where $\Phi_{0}$ is the vacuum vector and we use the convention that for any linear
operator $Y$, $Y^0=id$.\\
For simplicity of notations we \textbf{normalize the semi--circle principal Jacobi sequence}
\eqref{om01} so that
\begin{equation}\label{normaliz-om-n}
\omega = \omega_{n} =1 \qquad,\qquad\forall n\geq 1
\end{equation}
(see discussion in section \ref{sec:The-mu-Hilb-transf}).
With this normalization $\|\Phi_{n}\|=\omega_{n}!=1$, i.e. \textbf{the monic polynomials coincide with the normalized polynomials} and are an ortho--normal basis of $\Gamma_{X}$.
Recall that we use the convention
\begin{equation}\label{noMoSc01c0}
\Phi_{-n}=0 \qquad; \quad \ \forall n\in\mathbb{N}^*
\end{equation}
and, because of \eqref{normaliz-om-n},
\begin{equation}\label{act-a-a+1}
a\Phi_{n}=\Phi_{n-1} \quad;\quad a^{+}\Phi_{n} =\Phi_{n+1} \quad;\quad
[a,a^{+}] = \partial\omega_{\Lambda}
\end{equation}
The following Lemma recalls some properties of the canonical quantum decomposition of the
semi--circle random variable.
\begin{lemma}{\rm
In the canonical representation of the semi--circle random variable
with $\omega=1$, the following multiplication table holds:
\begin{equation}\label{free-MT}
\omega_{\Lambda+1} = aa^{+} = 1   \quad;\quad
\omega_{\Lambda}\Big|_{\{\Phi_{0}\}^{\perp}} =  a^{+}a\Big|_{\{\Phi_{0}\}^{\perp}} = 1
\quad;\quad a^{+}a\Phi_{0} = 0
\end{equation}
In particular
\begin{equation}\label{d-omL-SC}
\omega_{\Lambda} = a^{+}a = 1-\Phi_{0}\Phi_{0}^*
\end{equation}
where, for any $\xi\in\Gamma_{X}$, $\xi^*$ denotes the linear functional
$\xi^*:\eta\in\Gamma_{X}\to \xi^*(\eta) := \langle\xi, \eta\rangle$,
\begin{equation}\label{d-om-Lam-SC}
\partial\omega_{\Lambda} = \Phi_{0} \Phi_{0}^* = \delta_{0,\Lambda}
\end{equation}
\begin{equation}\label{free-CR}
[a,a^{+}] = \partial\omega_{\Lambda}= \Phi_{0} \Phi_{0}^*
\end{equation}
}\end{lemma}
\textbf{Proof}.
The first identity in \eqref{free-MT} follows from \eqref{act-a-a+1} because
$$
aa^{+}\Phi_{n} = a\Phi_{n+1} = \Phi_{n}
$$
The third identity in \eqref{free-MT} follows from the Fock property, i.e.
$a\Phi_{0} = 0$. For $n\ge 1$ one has
$$
a^{+}a\Phi_{n} = a^{+}\Phi_{n-1} = \Phi_{n}
$$
Thus $a^{+}a\Big|_{\{\Phi_{0}\}^{\perp}} = 1 $ which is the second identity in \eqref{free-MT}.
The last two identities in \eqref{free-MT} imply \eqref{d-omL-SC} and, given this
$$
[a,a^{+}]
= 1 - (1-\Phi_{0}\Phi_{0}^*)
= \Phi_{0}\Phi_{0}^*
$$
which is \eqref{free-CR}.
This implies the first identity in \eqref{d-om-Lam-SC} because of \eqref{comm-rel-[a,a+]}.\\
The second identity in \eqref{d-om-Lam-SC} follows from the definition of $F_{\Lambda}$ with $F_{\Lambda}=\delta_{0,\Lambda}$.
$\qquad\square$\\

\noindent\textbf{Remark}.
Note that \eqref{d-om-Lam-SC} implies that, for  $m\ge 1$
$$
[a,a^{+}]\Phi_{n} = \partial\omega_{\Lambda}\Phi_{n} = 0
\quad,\quad\forall n\ge 1
$$
Thus, in the canonical representation of the semi--circle law,
\textbf{the\\ non--commutativity of $a$ and $a^{+}$ is restricted to the vacuum space}.\\

\subsection{The $*$--Lie--algebra associated to the standard semi--circle distribution}

In this section we prove an infinite dimensional extension of Theorem 8 in \cite{[AcBaLuRha15]}.
The proof is a drastic simplification of the computational proof given in that paper.\\
Denote
\begin{equation}\label{df-P-Phin}
P_{\Phi_{m}}:= \hbox{ the projection onto } \mathbb{C}\cdot\Phi_{m}
\  \ \left(P_{\Phi_{m}}:=\delta_{m,\Lambda}=\Phi_{m}\Phi_{m}^*\right)
\end{equation}
With this notation one has
\begin{equation}\label{om05a}
aP_{\Phi_{0}} = P_{\Phi_{0}}a^{+}=0
\end{equation}
and the associated commutation relations are
\begin{equation}\label{om05a2}
 \left[ a,a^{+}\right] = \partial\omega_{\Lambda} = P_{\Phi_{0}}
\quad;\quad \left[  a,P_{\Phi_{0}}\right]
=-P_{\Phi_{0}}a\quad;\quad \ \ \left[  a^{+},P_{\Phi_{0}}\right]  =a^{+}P_{\Phi_{0}}
\end{equation}
\begin{theorem}\label{om05}{\rm
The $*$--Lie algebra generated by $a$ and $a^{+}$ is:
\begin{equation}\label{L0-SC}
\mathcal{L}_{0} = (\mathbb{C}\cdot a)\oplus (\mathbb{C}\cdot a^{+})\oplus
\mathcal{L}_{rank \, 1}\left(\Gamma_{X}\right)
\end{equation}
where $\mathcal{L}_{rank \, 1}\left(\Gamma_{X}\right)$ denotes the
$*$ algebra of rank $1$ operators on $\Gamma_{X}$ generated by the $(\Phi_{n})$.
}\end{theorem}
\textbf{Proof}.
By definition $a,a^{+} \in\mathcal{L}_{0}$ and, by \eqref{om05a2}, $P_{\Phi_{0}}\in\mathcal{L}_{0}$
and, $P_{\Phi_{0}}a\in\mathcal{L}_{0}$. Suppose by induction that
$P_{\Phi_{0}}a^{m}\in\mathcal{L}_{0}$. Then
\begin{equation}\label{om05a2a}
[a,P_{\Phi_{0}}a^{m}]
=[a,P_{\Phi_{0}}]a^{m}
\overset{\eqref{om05a2}}{=} \ - P_{\Phi_{0}}a^{m+1}\in\mathcal{L}_{0}.
\end{equation}
Therefore by induction $P_{\Phi_{0}}a^{m}\in\mathcal{L}_{0}$ for each $m\in\mathbb{N}$.
Taking adjoints of \eqref{om05a2a}, one finds
\begin{equation}\label{om05a2b}
[a^{+m}P_{\Phi_{0}},a^{+}]
= - a^{+(m+1)}P_{\Phi_{0}}\in\mathcal{L}_{0}.
\end{equation}
Therefore
$$
\mathcal{L}_{0}\ni
[a^{+m}P_{\Phi_{0}}, P_{\Phi_{0}}a^{n}]
=[a^{+m}P_{\Phi_{0}}, P_{\Phi_{0}}]a^{n}
+P_{\Phi_{0}}[a^{+m}P_{\Phi_{0}}, a^{n}]
$$
$$
=[a^{+m}, P_{\Phi_{0}}]P_{\Phi_{0}}a^{n}
+P_{\Phi_{0}}a^{+m}[P_{\Phi_{0}}, a^{n}]
+P_{\Phi_{0}}[a^{+m}, a^{n}]P_{\Phi_{0}}
$$
$$
\overset{\eqref{om05a2a}, \eqref{om05a2b}}{=} \
a^{+m}P_{\Phi_{0}}a^{n}
+P_{\Phi_{0}}a^{+m}P_{\Phi_{0}}a^{n}
+P_{\Phi_{0}}[a^{+m}, a^{n}]P_{\Phi_{0}}
$$
$$
\overset{\eqref{df-P-Phin}}{=}
a^{+m}P_{\Phi_{0}}a^{n}
+\langle \Phi,a^{+m}\Phi\rangle P_{\Phi_{0}}a^{n}
+\langle \Phi, [a^{+m}, a^{n}] \Phi\rangle P_{\Phi_{0}}
$$
$$
= a^{+m}P_{\Phi_{0}}a^{n}
=\Phi_{m}\Phi_{n}^*.
$$
The linear space generated by the set $\{\Phi_{m}\Phi_{n}^* \ : \ m,n\in\mathbb{N}\}$ consists
of all finite rank operators, is closed under commutators, invariant under commutators by both
$a^{+}$ and $a$ and has zero intersection with both $\mathbb{C}\cdot a^{+}$ and
$\mathbb{C}\cdot a$.
It follows that the right hand side of \eqref{L0-SC} is a $*$--Lie algebra with linear generators
$a^{+}, a$ and $\{\Phi_{m}\Phi_{n}^* \ : \ m,n\in\mathbb{N}\}$.
On the other hand the above arguments show that any $*$--Lie algebra containing $a^{+}$ and $a$
must contain all these operators. This proves \eqref{L0-SC}.
$\qquad\square$

\section{Analytic forms of free momentum and free kinetic energy operator}\label{analytic-forms}

\subsection{The $\mu$--Hilbert transform}\label{sec:The-mu-Hilb-transf}

Let $\mu_a$ be the semi-circle measure supported on $[-a,a]$ and $\lambda$ be the Lebesgue measure
on $\mathbb{R}$. Recall that the Hilbert transform of a function on the interval
$I_a=[-a,a]\subseteq\mathbb{R}$ is defined by
\begin{equation}\label{df-Hilb-Tr}
H_{I_a}\left[\chi_{I_a}(y)f(y)\right](x):= \frac{1}{\pi}\hbox{p.v.}\int_{I_a}\frac{f(y)}{x-y} dy
\end{equation}
where $\chi_{I_a}(x)=1$ if $x\in I_a$, $=0$ if $x\notin I_a$.
The change of integration variable $x\mapsto x/a$ (corresponding to the re--scaling
$\omega\to \omega/a^{2}$ of the principal Jacobi sequence \eqref{om01}, see \cite{[HoOb06]} Proposition (1.49)) leads to the study of  $H_{I_1}$.
Recall from \cite[Proposition 8.1.9]{[BuNe71]} that $H_{I_1}$ is a bounded operator on $L^p(I_1,\lambda)$, for $1< p<\infty$ into itself.
Moreover, by [Corollary 2.3] in \cite{[APS96]}, it extends to a bounded \textbf{surjection}
on $L^2(\mathbb{R},{\mu_1})$ \textbf{with a dense image}.
We mention here that $L^2(I_1,\lambda)$ is a proper sub--space of $L^2(\mathbb{R},{\mu_1})$ since $1/\sqrt{1-t^2}\in L^2(\mathbb{R},{\mu_1})$ but is not in $L^2(I_1,\lambda)$.\\

\noindent In this paper we are interested in the weighed Hilbert transform over the interval
$I_2=[-2,2]$:
\begin{equation*}
H_{\mu}f(x):= H_{I_2} \left[ \chi_{I_2}(y)f(y)\sqrt{4-y^2}\right](x)
\end{equation*}
i.e. the Hilbert transform with respect to the semi--circle measure $\mu$ over the interval $I_2$.
Using (4.176) in \cite{[King09]} together with the fact that
$$
f\in L^2(\mathbb{R},\mu_2)
\Rightarrow f(t)\sqrt{4-t^2}\in L^2(I_2,\lambda)
$$
since
\begin{equation*}
	\int_{I_2}|f(t)|^2(4-t^2)dt \le 2\int_{I_2}|f(t)|^2\sqrt{4-t^2}dt<\infty
\end{equation*}
it follows that, $ H_{\mu}$ is \textbf{skew-adjoint} on $L^2([-2,2],\mu)$.\\


\subsection{Representation on  $L^2([-2,2],{\mu})$}

Since $\mu$ has bounded support,
the polynomials $(\Phi_{n})_{n\ge0}$ form a \textbf{complete orthogonal system} in $L^2([-2,2],{\mu})$
(see, e.g., \cite{[MaHa03]}).
Assuming that  $(\omega_{n})_{n\ge1}$ and $(\alpha_{n})_{n\ge0}$ are respectively the constants $1$
and $0$, the polynomials $(\Phi_{n})_{n\ge0}$ are explicitly given by:
\begin{equation}\label{orth-pol-SC-Cheb}
	\Phi_{n}(x)=\frac{\sin((n+1)\arccos(x/2))}{\sin(\arccos(x/2))}, \quad n\ge0
\end{equation}
and satisfy the following
\textbf{monic Jacobi relation}
\begin{equation}\label{monic-recur}
	x \Phi_{n}(x)=\Phi_{n+1}(x)+\Phi_{n-1}(x), \quad \forall x\in [-2,2]
\end{equation}
The monic Chebyshev polynomial of first kind $T_n$ given by
\begin{equation}\label{df-Tn}
T_n(x)=2\cos(n\arccos(x/2))
\end{equation}
 are related with the $\Phi_{n}$ through the following relations
\begin{equation}\label{connection2}
	T_{n+1}(x)=\Phi_{n+1}(x)-\Phi_{n-1}(x)
\end{equation}
\begin{equation}\label{connection3}
	2T_{n+1}(x)=xT_{n}(x)-\left(4-x^2\right)\Phi_{n-1}(x)
\end{equation}
\begin{equation}\label{connection4}
	T_{n+1}(x)=2\Phi_{n+1}(x)-x\Phi_{n}(x)
\end{equation}
and the two classes of polynomials are connected via the $\mu$--Hilbert transform $ H_{\mu}$
as follows (see \cite{[Tricomi57]})
\begin{equation}\label{connection}
	H_{\mu}\Phi_{n}=T_{n+1} \quad,\quad n\ge0.
\end{equation}
\begin{proposition}\label{freemom}{\rm
		The operators $P_X$ and $i H_{\mu}$ coincide on $L^2([-2,2],\mu)$, i.e. for any $f\in L^2([-2,2],\mu)$, one has
		\begin{equation}\label{free-momentum}
			P_{X}f(x)=i H_{\mu} f(x)
			= 2i\hbox{p.v.}\int_{-2}^2\frac{f(y)}{x-y} \mu(dy).
		\end{equation}
	 In particular, the free kinetic energy operator is given by
		\begin{equation}\label{free-kin-energ}
			E_{X}:=\frac{1}{2}P_{X}^2 = -\frac{1}{2} H_{\mu}^2
		\end{equation}
}\end{proposition}
\textbf{Proof}.
By uniqueness of the continuous extension, it suffices to show the equality \eqref{free-momentum}
for the $\Phi_{n}$'s. From the equalities
\begin{equation*}
	a^+\Phi_{n}=\Phi_{n+1}\quad, \quad a\Phi_{n}=\omega_{n}\Phi_{n-1}=\Phi_{n-1}
\end{equation*}
we deduce that
\begin{align*}
	P_{X}\Phi_{n}=i(a^+-a)\Phi_{n}&=i(\Phi_{n+1}-\Phi_{n-1})
	\overset{\eqref{connection2}}{=}
	i T_{n+1}(x)
	\overset{\eqref{connection}}{=}
	i H_{\mu}\Phi_{n}
\end{align*}
This proves \eqref{free-momentum}. Given \eqref{free-momentum}, \eqref{free-kin-energ} is clear.
$\qquad\square$

\subsection{Generalized Schr\"odinger representation}

In this section, we provide the analytical forms of the free position, free momentum and
free CAP operators in the generalized Schr\"odinger representation.
To that goal, we let $Q$ denote the operator of multiplication by the coordinate in $L^2([-2,2],\mu)$ and $V$ the isometry from $L^2([-2,2],\mu)$ into $L^2(\mathbb{R},\lambda)$ given by
\begin{equation*}
	f\in L^2([-2,2],\mu) \mapsto f(Q) \rho\in L^2(\mathbb{R},\lambda)
\end{equation*}
where $\rho\in L^2(\mathbb{R},\lambda)$ is defined by
\begin{equation*}
	\rho (x):= \frac{1}{\sqrt{2\pi}}(4-x^2)^{1/4} \chi_{[-2,2]}(x), \quad x\in\mathbb{R}
\end{equation*}
Notice that, with $1$ denoting the constant function $=1$,
\begin{equation}\label{VPhi0=rho}
V\Phi_{0} =V1=\rho \ \in L^2(\mathbb{R},\lambda)
\end{equation}
Define $a^{\varepsilon}$ ($\varepsilon \in \{+,0,-\}$) the free CAP operators and set
\begin{equation*}
	A^{\varepsilon}:=Va^{\varepsilon}V^*
\end{equation*}
Since $V$ is isometric, the $A^{\varepsilon}$ satisfy the free
commutation relations:
\begin{equation}\label{free-comm-rel}
	[A^{-},A^{+}]=V\partial\omega_{\Lambda}V^*
\overset{\eqref{d-om-Lam-SC}}{=} V\Phi_{0} \Phi_{0}^* V^*=(V\Phi_{0})(V\Phi_{0})^*=\rho\rho^*
\end{equation}
\begin{theorem}\label{Gen-Schr-repr-Semicir-clas}{\rm
		\begin{enumerate}
			\item The free position operator
			$X$ is mapped into the usual position operator $Q= VXV^*$ in $L^2(\mathbb{R},\lambda)$:
			\begin{equation}\label{us-pos-op1}
				Qf(x) =xf(x)\ \quad, \quad x\in\mathbb{R}
			\end{equation}
			\item The free momentum operator $P_{X}$ is mapped into the operator $P = VP_{X}V^*$
			in $L^2(\mathbb{R},\lambda)$ given by:
\begin{equation}\label{us-mom-op1}
				P =i\rho H_{\mu} \rho^{-1}
			\end{equation}
			where, here and in the following, $\rho^{-1}:=1/\rho$ is understood on the support
of $\rho$ and we use the same symbol for $\rho$ and the multiplication operator by $\rho$.
			\item The free CAP operators $a^{\varepsilon}$ ($\varepsilon \in \{+,0,-\}$) are mapped by $V ( \, \cdot \, )V^*$ into the following CAP operators
			\begin{align}\label{mu-CAP-ops}
				A^{+} = \frac{1}{2}\big(X +\rho H_{\mu}  \rho^{-1}\big)\\
				A^{-}  = \frac{1}{2}\big(X  -\rho H_{\mu}  \rho^{-1}\big)
			\end{align}
\item The free kinetic energy operator $E_{X}$ is mapped into the operator\\ $E := VE_{X}V^*$
			in $L^2(\mathbb{R},\lambda)$ given by:
			\begin{equation}\label{us-kin-energ--op1}
				E = \frac{1}{2}P^2 =-\frac{1}{2}\rho H_{\mu}^2 \rho^{-1}
			\end{equation}
		\end{enumerate}
}\end{theorem}
\textbf{Proof}.
By direct computations one has, for each $f\in L^2(\mathbb{R},\lambda)$
$$
Qf(x) = VXV^*f(x)
=VX [\rho^{-1}f](x)
=xf(x)
$$
which is \eqref{us-pos-op1}. Similarly
$$
Pf = VP_{X}V^*f
=VP_{X} [\rho^{-1}f]
=iV H_{\mu} [\rho^{-1}f ]
=i\rho H_{\mu} [\rho^{-1}f]
$$
which is \eqref{us-mom-op1}.
Next, using the equalities
\begin{equation}\label{df-P-Q}
	Q = A^{+}  + A^{-} \qquad;\qquad
	P:=i(A^{+} - A^{-})
\end{equation}
one gets
\begin{equation*}
	A^{+}=\frac{1}{2}(Q -iP)
	\qquad;\qquad
	A^{-}=\frac{1}{2}(Q +iP)
\end{equation*}
which, given \eqref{us-pos-op1} and \eqref{us-mom-op1}, is \eqref{us-kin-energ--op1}.
Finally \eqref{us-kin-energ--op1} follows from,
\begin{align*}
	P^2f &=iP\rho H_{\mu} [\rho^{-1}f]=-\rho H_{\mu}^2 [\rho^{-1}f ]
\end{align*}
\begin{corollary}{\rm
		In Schr\"odinger representation, the free commutation relations \eqref{free-comm-rel} become:
		\begin{equation}\label{us-CCR-PQ1}
			[Q , P]
			=2i\rho\rho^*  \overset{\eqref{[PX,X]}}{=} V[X,P_X]V^*
		\end{equation}
		In particular,
		\begin{equation*}
			[Q , P] (\Phi_{n}\rho)=
			\begin{cases}
				2\rho\ ; n=0\\
				0\ \quad; n\ge1
			\end{cases}
		\end{equation*}
}\end{corollary}
\textbf{Proof}.
This follows from
$$
[Q , P] \overset{\eqref{df-P-Q}}{=}[A^{+}  + A^{-} , i(A^{+} - A^{-})]
=-i[A^{+}, A^{-})] + i[A^{-} , (A^{+}]
\overset{\eqref{free-comm-rel}}{=} 2i \rho\rho^*
$$
An analytic proof is the following. For any $f\in L^2(\mathbb{R},\lambda)$, one has
$$
P Q  f(x)
=i\rho H_{\mu} [Qf/\rho](x)
=2i\rho \hbox{p.v.}\int_{-2}^2\frac{yf(y)}{x-y}\frac{1}{\rho(y)} \mu(dy)
$$
$$
=-2i\rho \hbox{p.v.}\int_{-2}^2\frac{(x-y)f(y)}{x-y}\frac{1}{\rho(y)} \mu(dy)
+2xi\rho \hbox{p.v.}\int_{-2}^2\frac{f(y)}{x-y}\frac{1}{\rho(y)} \mu(dy)
$$
$$
=-2i\rho \int_{-2}^2 f(y) \frac{1}{\rho(y)}\rho(y)^2 dy
+xi\rho H_{\mu}[\rho^{-1}f](x)
$$
$$
=-2i\rho \langle f,\rho\rangle + [Qi\rho  H_{\mu}\rho^{-1}f](x)
=-2i\rho \langle f,\rho\rangle + QPf(x)
$$
which is equivalent to \eqref{us-CCR-PQ1}.
$\qquad\square$\\


\section{Free harmonic oscillators in  generalized quantum mechanics}

The results in this section are valid \textbf{for any} principal Jacobi sequence $(\omega_{n})$.\\
Define Hamiltonian of the \textbf{semi--circle harmonic oscillator} with frequency $\hat{\omega}$ with the
same formula used for the usual $\hat{\omega}$--harmonic oscillator, i.e.
\begin{equation}\label{Harm-Osc-om}
H_{\hat{\omega}} :=  \frac{1}{2}(p^2 + \hat{\omega}^2 q^2)
\end{equation}
In generalized quantum mechanics this becomes
$$
\frac{1}{2}(p^2 + \hat{\omega}^2 q^2)
=\frac{1}{2}(-(a^{+} - a)^2 + \hat{\omega}^2 (a^{+} + a)^2)
$$
$$
=\frac{1}{2}(-(a^{+2} - a^{+}a - aa^{+} + a^2)
+ \hat{\omega}^2 (a^{+2} + a^{+}a + aa^{+} + a^2) )
$$
$$
=\frac{1}{2}((\hat{\omega}^2-1)(a^{+2} + a^2) + (\hat{\omega}^2+1)(a^{+}a + aa^{+}) )
$$
\begin{equation}\label{Harm-Osc-gen}
\overset{\eqref{a+a-aa+}}{=}
\frac{1}{2}((\hat{\omega}^2-1)(a^{+2} + a^2)
+ (\hat{\omega}^2+1)(\omega_{\Lambda} + \omega_{\Lambda+1}) )
\end{equation}
Putting $\hat{\omega}^2=1$ in \eqref{Harm-Osc-om}, one finds
\begin{equation}\label{Harm-Osc-freq=1}
H_{1} =  \omega_{\Lambda} + \omega_{\Lambda+1}
\end{equation}
In the semi--circle case with $\omega_{n}=1$ for each $n$, recalling the notation
\eqref{df-P-Phin} (i.e. $P_{\Phi_{0}}:=\Phi_{0}\Phi_{0}^*$), for any $n\in\mathbb{N}$, this implies
$$
e^{itH_{1}}
=  e^{it(\omega_{\Lambda} + \omega_{\Lambda+1})}
=  e^{it(\omega_{\Lambda} + \omega_{\Lambda+1})}P_{\Phi_{0}}
+ e^{it(\omega_{\Lambda} + \omega_{\Lambda+1})}P_{\Phi_{0}}^{\perp}
$$
\begin{equation}\label{Harm-Osc-freq=1a}
=  e^{it\omega_{1}}P_{\Phi_{0}}
+ e^{it2\omega_{1}}P_{\Phi_{0}}^{\perp}.
\end{equation}
Therefore, \textbf{for any unit vector}
$\xi=P_{\Phi_{0}}\xi + P_{\Phi_{0}}^{\perp}\xi =:\xi_{0}+\xi^{\perp}\in\Gamma_{X}$, one has
$$
\langle\xi, e^{itH_{1}} \xi\rangle
=\langle\xi, e^{it\omega_{1}}P_{\Phi_{0}}\xi + e^{it2\omega_{1}}P_{\Phi_{0}}^{\perp}\xi\rangle
=\langle\xi_{0}+\xi^{\perp}, e^{it\omega_{1}}\xi_{0} + e^{it2\omega_{1}}\xi^{\perp}\rangle
$$
$$
=\|\xi_{0}\|^2e^{it\omega_{1}} + \|\xi^{\perp}\|^2e^{it2\omega_{1}}.
$$
Thus, denoting
$$
p_{\xi} := |\xi_{0}|^2 \in [0,1]  \qquad;\qquad 1-p_{\xi} := \sum_{n=1}^{\infty}|\xi_{n}|^2
$$
one obtains
$$
\langle\xi, e^{it(\omega_{\Lambda}+\omega_{\Lambda+1})} \xi\rangle
= p_{\xi}e^{it\omega_{1}} + (1-p_{\xi})e^{it2\omega_{1}}
$$
i.e. the characteristic function of the \textbf{Bernoulli random variable} with values
$\{\omega_{1},2\omega_{1}\}$ and distribution $(p_{\xi}, 1-p_{\xi})$.\\

\noindent\textbf{Remark}.
Notice that the structure of the solution becomes considerably more complex in the case
$\hat{\omega}^2\ne 1$ in \eqref{Harm-Osc-om}. In this case in fact, the right hand side of
\eqref{Harm-Osc-gen} differs only by the constants in the linear combination from the expression
of the kinetic energy
$$
\frac{1}{2}p^{2}
=\frac{1}{2}\left( a^{+2} + a^2 \right)
+\frac{1}{2}(\omega_{\Lambda} + \omega_{\Lambda+1}) )
$$
and therefore one can expect that it can be dealt with using techniques similar to those used for
the kinetic energy operator (see sect. \eqref{sec:evol-eitP2}). However this problem will not
be considered in the present paper.\\

\noindent\textbf{Remark}.
The associated \textbf{coherent vectors}, in the semi--circle case, are given by
$$
\psi_{z} := \sum_{n= 0}^{\infty}\frac{z^{n}}{\sqrt{\omega_{n}!}}\Phi_{n}
= \sum_{n= 0}^{\infty} z^{n}\Phi_{n} \qquad,\qquad |z|< 1
$$
so that the associated kernel is given by the geometric series.
$$
\langle\psi_{u} , \psi_{v}\rangle
= \sum_{n= 0}^{\infty} (\bar{u}v)^{n}
=\frac{1}{1-\bar{u}v}
 \qquad,\qquad |u|, |v|< 1
$$

\section{Evolutions generated by the free momentum}

\subsection{Action of the Heisenbeg evolution generated by $P_{X}$ on the CAP operators}
\label{sec:Act-Heis-ev-gen-P-CAP-ops}

Since $X$ is fixed, from now on we write simply
$$
P:=i (a^{+}-a)
$$
instead of $P_X$. Since $X$ is bounded, $P$ which is unitarily isomorphic to $X$, is bounded, hence
the $1$--parameter unitary group $e^{itP}$ generated by $P$ is well defined and its
exponential series converges in norm.
\begin{lemma}{\rm
		\begin{equation}\label{a+t-SC}
			a^{+}_t :=
			e^{itP}a^{+}e^{-itP}
			= a^{+} + \int_{0}^{t}ds (e^{isP_{X}}\Phi_{0})(e^{isP_{X}}\Phi_{0})^*
		\end{equation}
}\end{lemma}
\textbf{Proof}.
$$
\partial_ta^{+}_t
=ie^{itP_X}[P_X,a^{+}]e^{-itP_X}
= - e^{itP_X}[a^{+}-a,a^{+}]e^{-itP_X}
$$
$$
=e^{itP_X}i(-i\partial\omega_{\Lambda})e^{-itP_X}
=:\partial\omega_{\Lambda,t}
$$
Therefore, if $\partial\omega_{\Lambda}$ has the form \eqref{d-omL-SC},
$$
\partial_ta^{+}_t
= (\partial\omega_{\Lambda})_{t}
= (\Phi_{0} \Phi_{0}^*)_{t}
= e^{itP_{X}}(\Phi_{0} \Phi_{0}^*)e^{-itP_{X}}
$$
\begin{equation}\label{partial-t-a+t-SC}
	= (e^{itP_{X}}\Phi_{0})(e^{itP_{X}}\Phi_{0})^*
	\qquad;\qquad a^{+}_0 =a^{+}
\end{equation}
This is equivalent to
$$
a^{+}_t = a^{+} + \int_{0}^{t}ds (e^{isP_{X}}\Phi_{0})(e^{isP_{X}}\Phi_{0})^*
$$
which is \eqref{a+t-SC}.
$\qquad\square$\\

\textbf{Remark}.
For $n\in\mathbb{N}$,
\begin{equation}\label{a+t-SC2}
a^{+}_t\Phi_{n} \
\overset{\eqref{a+t-SC}}{=}  \
\Phi_{n+1} + \int_{0}^{t}ds \langle e^{isP_{X}}\Phi_{0} ,\Phi_{n}\rangle e^{isP_{X}}\Phi_{0}
\end{equation}
Therefore $a^{+}_t$, hence the action of $e^{isP_{X}} ( \, \cdot \,) e^{-isP_{X}}$ on the
$*$--algebra generated by $a^{+}$ and $a$, is completely determined by $e^{isP_{X}}\Phi_{0}$.

\subsection{Action of the Schr\"odinger evolution generated by $P_{X}$ on the number vectors}
\label{sec:Act-Schr-ev-gen-P-Nbr-vects}

We want to compute
$$
e^{itP}\Phi_{n}
=\sum_{k\ge 0}\frac{(it)^{k}}{k!}P^{k}\Phi_{n}
=\sum_{k\ge 0}\frac{(it)^{k}}{k!}(i)^{k}(a^{+}-a)^{k}\Phi_{n}
$$
\begin{equation}\label{eitP-Phi-n}
	=\sum_{k\ge 0}\frac{(-t)^{k}}{k!} \sum_{\varepsilon\in\left\{+,-\right\}^{k}}
	(-1)^{|\{j:\varepsilon(j)=-1\}|}
	a^{\varepsilon(k)}\cdots a^{\varepsilon(1)}\Phi_{n}
	\ ;\qquad \forall n\in\mathbb{N}
\end{equation}
Thus the problem is to evaluate the products
\begin{equation}\label{a-eps-Phi-n}
a^{\varepsilon(k)}\cdots a^{\varepsilon(1)}\Phi_{n}
\end{equation}
A standard way to attack this problem is to reduce the products of the form
\begin{equation}\label{a-eps}
a^{\varepsilon(k)}\cdots a^{\varepsilon(1)}
\end{equation}
to their normally ordered form, i.e. with creators on the left and annihilators on the right.
By iterated use of the commutation relation $aa^{+} = 1 $ (see \eqref{free-MT}), any
product \eqref{a-eps} can be reduced to the form
\begin{equation}\label{free-NO1}
a^{\varepsilon(k)}\cdots a^{\varepsilon(1)}
=(a^{+})^{m_{+}} a^{m_{-}}
\end{equation}
where the numbers $m_{+}, m_{-}\in\mathbb{N}$ are uniquely determined by $\varepsilon$ and,
when $\varepsilon$ varies in $\left\{+,-\right\}^{k}$, are all possible pairs satisfying
\begin{equation}\label{m+-m-}
m_{+} + m_{-} \in \{0,1,\dots,k\}
\end{equation}
So we know that $e^{itP}$ can be written in the form
\begin{equation}\label{rep1-PQ-mom-01d0}
e^{itP}
=\sum_{m,n\ge 0}I_{m,n}(t)(a^+)^{m}(a^-)^n
\end{equation}
for some numerical coefficients $I_{m,n}(t)$.
The problem is to calculate these coefficients. The solution of this problem requires the solution of
the \textbf{inverse normal order problem}, namely:
\textit{parametrize the set of $\varepsilon\in\left\{+,-\right\}^{k}$ that satisfy the
identity \eqref{free-NO1} when $k$ varies in $\mathbb{N}$} (c.f. \cite{[AcHamLu21]}).\\

\section{Free momentum group and free translations}\label{sec:evol-eitP}

In this section we discuss an application of Theorem 3.2 and Corollary 3.3 of Part II (\cite{[AcHamLu21]})
to the $1$MIFS $\Gamma\big(\mathbb{C},\{\omega_{n}\}_{n\ge 1}\big)$
with $\omega_{n}=1$ for any $n\ge 1$.\\
In this case, one knows that the spectrum of the \textit{position operator }$X:=a^-+a^+$ is
the interval $[-2,2]$. In particular $X$ is bounded and
\begin{equation}\label{norm-X-SC}
\|X\| = 2
\end{equation}

\subsection{The vacuum distribution of $P$ and $X$ }

From the general theory of orthogonal polynomials we know that
$$
P=\Gamma(i) X\Gamma(i)^*
$$
where $\Gamma(i)$ is the Gauss--Fourier transform associated to $X$, which is a unitary operator
(see \cite{[AcEllLu20-mom]}).
Therefore $X$ and $P=P_X$ have the same spectrum and, since
$\Gamma(i)\Phi_{0}=\Phi_{0}= \Gamma(i)^*\Phi_{0}$ they have the same vacuum distribution
(this is true for every classical random variable with all moments. The following Lemma
gives a direct combinatorial proof in case of the semi--circle law).
\begin{lemma}\label{X-P-same-vac-distr}{\rm
The vacuum distribution of $X$ and of the $X$--momentum operator
\begin{equation}\label{notat-for-P}
P:=-i(a^--a^+)
\end{equation}
coincide and are the semi--circle distribution on the interval $[-2,2]$.
}\end{lemma}
\textbf{Proof}.
It is sufficient to prove that
\begin{equation}\label{P-mom}
\langle \Phi_{0}, P^{2n} \Phi_{0}\rangle = C_{n} \hbox{ (the $n$--th Catalan number)}, \quad\forall n\in\mathbb{N}
\end{equation}
Denoting
\begin{equation}\label{df-nu-pm(eps)}
\nu_\pm(\varepsilon):=|\varepsilon^{-1}(\{\pm\})| \quad;\quad n\in\mathbb{N} \ , \
\varepsilon\in\{-,+\}^n
\end{equation}
one has
$$
\nu_\pm(\varepsilon)=m , \hbox{ if }\varepsilon\in\{-,+\}^{2m}_+
$$
\begin{equation}\label{PQ-mom-00}
\begin{cases}
X^n=\sum_{\varepsilon\in\{-,+\}^n}a^{\varepsilon(1)}\cdots a^{\varepsilon(n)} \\
P^n\overset{\eqref{notat-for-P}}{=}
\sum_{\varepsilon\in\{-,+\}^n} (-i)^n (-1)^{\nu_+(\varepsilon)}
a^{\varepsilon(1)}\cdots a^{\varepsilon(n)}
\end{cases}
\end{equation}
Denoting $\Phi_0$:=the vacuum vector, one has
\begin{equation}\label{PQ-mom-00a}
\langle \Phi_0,a^{\varepsilon(1)}\cdots a^{\varepsilon(n)}\Phi_0\rangle
=\begin{cases}
   1, & \mbox{if } n=2m \hbox{ and  }\varepsilon\in\{-,+\}^{2m}_+ \\
   0, & \mbox{otherwise}.
 \end{cases}
\end{equation}
Consequently
\begin{align}\label{PQ-mom-00b}
&\langle \Phi_0,P^n\Phi_0\rangle
=\sum_{\varepsilon\in\{-,+\}^n} (-i)^n (-1)^{\nu_+(\varepsilon)}\langle\Phi_0,a^{\varepsilon(1)}
\cdots a^{\varepsilon(n)}\Phi_0\rangle\\
&\overset{\eqref{PQ-mom-00a},\eqref{PQ-mom-00}}{=}
\begin{cases}
\sum_{\varepsilon\in\{-,+\}^{2m}_+}(i)^{2m} (-1)^m, & \mbox{if } n=2m \\
0, & \mbox{ if $n$ is odd }
\end{cases}
\notag\\
&=\begin{cases}
|\{-,+\}^{2m}_+|, & \mbox{if } n=2m \\
0, & \mbox{ if $n$ is odd }
  \end{cases}
 \ = \begin{cases}
C_m, & \mbox{if } n=2m \\
0, & \mbox{ if $n$ is odd }
     \end{cases} \\
&= \langle \Phi_0,X^n\Phi_0\rangle \notag
\end{align}
which are the moments of the semi--circle distribution on the interval $[-2,2]$.
Since for this distribution, the moment problem is determined, the two distributions coincide.
This proves \eqref{P-mom}.
$\qquad \square$\\

\subsection{The evolutions $e^{itP}$, $e^{itX}$}\label{sec:The-ev-eitP-eitX}

\noindent Since $X$ and $P=P_X$ have the same spectrum, \eqref{norm-X-SC} holds with
$X$ replaced by $P$ hence, for all $z\in\mathbb{C}$, the exponential series of both $e^{zP}$
and $e^{zX}$ converge in norm.\\
In this section we will determine the action of the groups $e^{zP}$ and $e^{zX}$ on the
monic basis \eqref{df-mon-bas}.
\begin{theorem}\label{PQ-mom-01}{\rm
In the notations of Theorem 3.2, Proposition 3.6, and Corollary 3.3 of Part II
(\cite{[AcHamLu21]}) and assuming that $\omega_{n}=1$ for any $n\ge 1$, one has,
\begin{equation}\label{df-Imn(t)}
I_{m,n}(t):=\sum_{p\ge 0}\frac{t^{m+n+2p}}{(m+n+2p)!}
(-1)^{p+m}|\Theta_{m+n+2p}(m,n)|
\end{equation}
where, for any $m_+, m_-, p\ge 0$,
\begin{equation}\label{|Theta(m++m-+2p)(m+,m-)|}
|\Theta_{m_++m_-+2p}(m_+,m_-)|=\frac{m_++m_-+1}{2p+m_++m_-+1}\binom{2p+m_++m_-+1}{p}
\end{equation}
and for any $t\in\mathbb{C}$,
\begin{equation}\label{rep1-PQ-mom-01d}
e^{itP}
=\sum_{m,n\ge 0}I_{m,n}(t)(a^+)^{m}(a^-)^n.
\end{equation}
}\end{theorem}
\textbf{Proof}.
In order to prove \eqref{rep1-PQ-mom-01d} notice that for any choice of $m,m_+,m_-$ with $m\ge m_++m_-$,
one has
\begin{align*}
\{-,+\}^{m}
& =
\bigcup_{m_+,m_-\in\mathbb{N}; m-m_+-m_-\hbox{ even}}\Theta_m(m_+,m_-)\\
&=\bigcup_{m_+,m_-,p\ge 0; m_++m_-+2p=m}\Theta_{m_++m_-+2p}(m_+,m_-)
\end{align*}
and the union is disjoint. Consequently
\begin{equation}\label{m-to-m+m-p}
\bigcup_{m\ge 0}\{-,+\}^{m}=\bigcup_{m_+,m_-,p\ge 0}\Theta_{m_++m_-+2p}(m_+,m_-)
\end{equation}
both unions being disjoint. This implies
$$
e^{itP} =\sum_{m\ge 0}\frac{(it)^{m}}{m!}P^{m}
\overset{\eqref{PQ-mom-00}}{=}
\sum_{m\ge 0}\sum_{\varepsilon\in\{-,+\}^{m}}\frac{(i)^{m}t^{m}}{m!}(-i)^{m}
(-1)^{\nu_+(\varepsilon)}a^{\varepsilon(1)}\cdots a^{\varepsilon( m)  }
$$
\begin{equation}\label{itP-2}
=\sum_{m\ge 0}\sum_{\varepsilon\in\{-,+\}^{m}}\frac{t^{m}}{m!}
(-1)^{\nu_+(\varepsilon)}a^{\varepsilon(1)}\cdots a^{\varepsilon( m)  }
\end{equation}
$$
\overset{\eqref{m-to-m+m-p}}{=}\sum_{m_+,m_-,p\ge 0}\frac{t^{m_++m_-+2p}}{(m_++m_-+2p)!}
\sum_{\varepsilon\in\Theta_{m_++m_-+2p}(m_+,m_-)}
(-1)^{\nu_+(\varepsilon)}a^{\varepsilon(1)}\cdots a^{\varepsilon(m_++m_-+2p)  }.
$$
From Corollary 3.3 of \cite{[AcHamLu21]}
we know that for any $m_+,m_-,p\ge 0$,
\begin{equation}\label{PQ-mom-01c}
a^{\varepsilon(1)}\cdots a^{\varepsilon( m_++m_-+2p)  }=(a^+)^{m_+}(a^-)^{m_-};\ \ \forall \varepsilon\in\Theta_{m_++m_-+2p}(m_+,m_-)
\end{equation}
and $\nu_+$ takes value $p+m_+$ on $\Theta_{m_++m_-+2p}(m_+,m_-)$. Therefore
\begin{align*}
e^{itP} &= \sum_{m_+,m_-,p\ge 0}
\frac{(-1)^{p+m_+}t^{m_++m_-+2p}}{(m_++m_-+2p)!}  
|\Theta_{m_++m_-+2p}(m_+,m_-)|(a^+)^{m_+}(a^-)^{m_-}
\label{PQ-mom-01g}
\\&= \sum_{m_+,m_-\ge 0}\big(\sum_{p\ge 0}
\frac{(-1)^{p+m_+}t^{m_++m_-+2p}}{(m_++m_-+2p)!}
|\Theta_{m_++m_-+2p}(m_+,m_-)|\big)(a^+)^{m_+}(a^-)^{m_-}.
\end{align*}
Using the definition of the $I_{m,n}(t)$ given by the identity  \eqref{df-Imn(t)}
this becomes
$$
e^{itP}
=\sum_{m,n\ge 0}I_{m,n}(t)(a^+)^{m}(a^-)^n
$$
which is the identity  \eqref{rep1-PQ-mom-01d}.
$\qquad\square$\\
\begin{lemma}\label{df-Imn(t)-01}
\noindent Let
\begin{equation}\label{df-Bessl-fctn-1rst-kind}
	J_n(t):=\sum_{p\ge0}\frac{(-1)^{p}}{p!(n+p)!}\left(\frac{t}{2}\right)^{n+2p}
\end{equation}
denote the Bessel function of first kind.
For any $m,n\in\mathbb{N}$, we have 
\begin{align*}
I_{m,n}(t)=(-1)^{m}I_{0,m+n}(t)=(-1)^{n}I_{m+n,0}(t)
\end{align*}
where
\begin{align}\label{I0n(t)}
	I_{0,n}(t)=(n+1)\frac{J_{n+1}(2t)}{t}=J_{n+2}(2t)+J_{n}(2t)
\end{align}
\end{lemma}
\textbf{Proof}.
From the definition of $I_{m,n}(t)$ given by the first identity in \eqref{df-Imn(t)}, it is
clear that for $t=0$,
\begin{equation*}
I_{m,n}(0)=\delta_{m+n,0}
\end{equation*}
Now, if $t\in\mathbb{C}\setminus\{0\}$, one has
\begin{align*}
I_{m,n}(t)&=\sum_{p\ge0}\frac{t^{2p+m+n}}{(2p+m+n)!}(-1)^{p+m}|\Theta_{m+n+2p}(m,n)|
\\&\overset{\eqref{|Theta(m++m-+2p)(m+,m-)|}}{=}
\sum_{p\ge0}\frac{(-1)^{p+m}t^{2p+m+n}}{(2p+m+n)!}\frac{m+n+1}{2p+m+n+1}\binom{2p+m+n+1}{p}
\\&=\sum_{p\ge0}\frac{(-1)^{p+m}t^{2p+m+n}}{p!}\frac{m+n+1}{(p+m+n+1)!}
\\&=(-1)^{m}(m+n+1)\frac{J_{m+n+1}(2t)}{t}.
\end{align*}
Then, using the recurrence relation
\begin{align*}
\frac{2nJ_{n}(t)}{t}=J_{n+1}(t)+J_{n-1}(t)
\end{align*}
we obtain the result.
$\qquad\square$
\begin{proposition}\label{exp-evol-P}
For all $t\in\mathbb{C},k\in\mathbb{N}$ and $x\in[-2,2]$,
\begin{equation}
e^{itP}\Phi_k(x)=J_{0}(2t)\Phi_{k}(x)
-\sum_{n=1}^{k}J_{n}(2t)T_{k-n+2}(x)
\end{equation}
$$
+x\sum_{n=0}^{k}\sum_{m\ge n+2}(-1)^{m-n}J_{m}(2t)\Phi_{m+k-2n-1}(x)
$$
\end{proposition}
\textbf{Proof}.
From Theorem \ref{PQ-mom-01} and Lemma \ref{df-Imn(t)-01}, we have
\begin{align*}
e^{itP}&=\sum_{m,n\ge0}(-1)^{m}(J_{m+n+2}(2t)+J_{m+n}(2t))(a^+)^m(a^-)^n
\end{align*}
Then, using the fact that for all $m,n,k\ge0$
\begin{equation*}
(a^+)^m(a^-)^n\Phi_k=
\begin{cases}
 \Phi_{m+k-n}, & \mbox{if} \ n\le k\\
 0 , & \mbox{otherwise}
\end{cases}
\end{equation*}
we obtain
\begin{align*}
e^{itP}\Phi_k(x)=&\sum_{n=0}^{k}\sum_{m\ge0}(-1)^{m}(J_{m+n+2}(2t)+J_{m+n}(2t))\Phi_{m+k-n}(x)
\\=&\sum_{n=0}^{k}\sum_{m\ge n}(-1)^{m-n}(J_{m+2}(2t)+J_{m}(2t))\Phi_{m+k-2n}(x)
\\=&\sum_{n=0}^{k}J_{n}(2t)\Phi_{k-n}(x)-J_{n+1}(2t)\Phi_{k+1-n}(x)
\\&+\sum_{m\ge n+2}(-1)^{m-n}J_{m}(2t)\left\{\Phi_{m+k-2n}(x)+\Phi_{m+k-2n-2}(x) \right\}
\\=&J_{0}(2t)\Phi_{k}(x)+\sum_{n=1}^{k}J_{n}(2t)[\Phi_{k-n}(x)-\Phi_{k+2-n}(x)]
\\&+x\sum_{n=0}^{k}\sum_{m\ge n+2}(-1)^{m-n}J_{m}(2t)\Phi_{m+k-2n-1}(x)
\\=&J_{0}(2t)\Phi_{k}(x)-\sum_{n=1}^{k}J_{n}(2t)T_{k-n+2}(x)
\\&+x\sum_{n=0}^{k}\sum_{m\ge n+2}(-1)^{m-n}J_{m}(2t)\Phi_{m+k-2n-1}(x)
\end{align*}
which proves the proposition. $\qquad\square$\\

\begin{proposition}
	For any $t\in\mathbb{C}$ and any $k\in\mathbb{N}$,
		\begin{align*}
		e^{itP}\Phi_k(x)=\sum_{l=0}^{\infty}\langle \Phi_l,e^{itP}\Phi_k\rangle\Phi_{l}(x)
	\end{align*}
where
	\begin{equation*}
	\langle \Phi_l,e^{itP}\Phi_k\rangle=\sum_{m=0}^{l\wedge k}I_{l-m,k-m}(t)
\end{equation*}
In particular,
\begin{equation}\label{eitP-phi0}
	e^{itP}\Phi_0(x)=\sum_{l=0}^{\infty}I_{l,0}(t)\Phi_{l}(x)=\frac{1}{t}
\sum_{l=0}^{\infty}(-1)^l(l+1)J_{l+1}(2t)\Phi_{l}(x)
\end{equation}
\end{proposition}
\textbf{Proof}.
	For any $t\in\mathbb{C}$, one has
	\begin{align*}
		e^{itP}\Phi_k(x)=&\sum_{m=0}^{\infty}\sum_{n=0}^{k}I_{m,n}(t)\Phi_{m+k-n}(x)
=\sum_{m=0}^{\infty}\sum_{n=0}^{k}I_{m,k-n}(t)\Phi_{m+n}(x)
	\end{align*}
	Performing the variable change $l=m+n$ and rearranging terms, one gets
	\begin{align*}
		e^{itP}\Phi_k(x)=&\sum_{l=0}^{\infty}\sum_{m=0}^{l}I_{m,k-(l-m)}(t)\chi_{\{k\ge l-m\}}\Phi_{l}(x)=\sum_{l=0}^{\infty}\sum_{m=0}^{l\wedge k}I_{l-m,k-m}(t)\Phi_{l}(x)
	\end{align*}
	Hence, for all $k,l\in\mathbb{N}$
	\begin{equation*}
		\langle \Phi_l,e^{itP}\Phi_k\rangle=\sum_{m=0}^{l\wedge k}I_{l-m,k-m}(t)
	\end{equation*}
	In particular,
	\begin{equation*}
		\langle \Phi_l,e^{itP}\Phi_0\rangle=I_{l,0}(t)
=(-1)^l(l+1)\frac{J_{l+1}(2t)}{t}=(-1)^l\langle \Phi_0,e^{itP}\Phi_l\rangle
	\end{equation*}
which proves \eqref{eitP-phi0}. $\qquad\square$\\

\begin{remark}

Notice that the \textbf{characteristic function} of $P$ with respect to the state $\langle \Phi_l, \ \cdot \ \Phi_l\rangle$ on $\mathcal{B}(\Gamma\big(\mathbb{C},\{\omega_{n}\}_{n\ge 1}))$ is
\begin{equation*}
\sum_{m=0}^{l}I_{m,m}(t)
\end{equation*}
In particular, the vacuum expectations of $e^{itP}$ is given by
\begin{equation}\label{vac-exp-eitP}
\langle\Phi_0,e^{itP}\Phi_0\rangle=I_{0,0}(t)
=\sum_{p=0}^\infty\frac{(-t^2)^{p}}{(2p)!} C_p=\frac{J_1(2t)}{t}
\end{equation}
\end{remark}
\begin{lemma}\label{Bessel-sum}
We have
\begin{equation*}
\sum_{m\ge 1}(-1)^mJ_{m}(2t)\sin(m\theta)
=-\frac{\sin(2t\sin\theta)}{2}
-\frac{\sin(\theta)}{2\pi}\int_0^{\pi}\frac{\cos(2t\sin\varphi)}{\cos(\theta)-\cos(\varphi)}d\varphi
\end{equation*}
and
\begin{align*}
\sum_{m\ge 1}(-1)^mJ_{m}(2t)\cos(m\theta)=&\frac{\cos(2t\sin\theta)-J_{0}(2t)}{2}
-\frac{1}{2\pi}\int_0^{\pi}\frac{\sin(2t\sin\varphi)\sin(\varphi)}{\cos(\theta)-\cos(\varphi)}d\varphi
\end{align*}
\end{lemma}
\textbf{Proof}.
Since these series are absolutely convergent, we can split them into sums of odd and even terms.
On the one hand, we have from \cite{[Kap1905]},
 \begin{equation*}
\sum_{m\ge 1}J_{2m}(2t)\sin(2m\theta)=-\frac{\sin(\theta)}{4\pi}\int_0^{2\pi}\frac{\cos(2t\sin\varphi)
-\cos(2t\sin\theta)}{\cos(\varphi)+\cos(\theta)}d\varphi
\end{equation*}
and
\begin{equation*}
\sum_{m\ge 0}J_{2m+1}(2t)\cos((2m+1)\theta)=\frac{1}{4\pi}\int_0^{2\pi}\frac{\sin(2t\sin\varphi)\sin(\varphi)
-\sin(2t\sin\theta)\sin(\theta)}{\cos(\varphi)+\cos(\theta)}d\varphi
\end{equation*}
Now, using \cite[(11.17)]{[King09]} together with the fact that the integrands are periodic and even functions, one obtains
\begin{equation*}
	\sum_{m\ge }J_{2m}(2t)\sin(2m\theta)=-\frac{\sin(\theta)}{2\pi}\int_0^{\pi}\frac{\cos(2t\sin\varphi)}{\cos(\theta)
-\cos(\varphi)}d\varphi
\end{equation*}
and
\begin{equation*}
	\sum_{m\ge 0}J_{2m+1}(2t)\cos((2m+1)\theta)=\frac{1}{2\pi}\int_0^{\pi}
\frac{\sin(2t\sin\varphi)\sin(\varphi)}{\cos(\theta)-\cos(\varphi)}d\varphi
\end{equation*}
On the other hand, by replacing $\theta$ with $\pi/2-\theta$, the Jacobi-Anger identity (\cite{CPVWJ08} p. 344):
\begin{align}\label{Jacobi-Anger}
e^{2it\cos\theta}
=\sum_{m=-\infty}^\infty i^mJ_{m}(2t)e^{im\theta}
=J_{0}(2t)+2\sum_{m\ge 1}i^{m}J_{m}(2t)\cos(m\theta)
\end{align}
becomes
\begin{align*}
e^{2it\sin\theta}=\sum_{m=-\infty}^\infty J_{m}(2t)e^{im\theta}
\end{align*}
Then, equating the real and imaginary parts respectively, we obtain
\begin{equation*}
\cos(2t\sin\theta)=J_{0}(2t)+2\sum_{m\ge 1}J_{2m}(2t)\cos(2m\theta)
\end{equation*}
and
\begin{equation*}
\sin(2t\sin\theta)=2\sum_{m\ge 0}J_{2m+1}(2t)\sin((2m+1)\theta)
\end{equation*}

which proves the Lemma. $\qquad\square$\\
\begin{proposition}
Let $x\in[-2,2]$ and $\theta=\arccos(x/2)$. Then, one has, for any $t\in\mathbb{C}$,
\begin{equation*}
		e^{itP}\phi_0(x)=J_0(2t)+xJ_1(2t)-\frac{x \sin(2t\sin(\theta))}{2\sin(\theta)}-\frac{x}{2\pi}\int_0^\pi\frac{\cos(2t\sin\varphi)}{\cos(\theta)
-\cos(\varphi)}d\varphi
\end{equation*}
and
\begin{equation*}
	e^{itP}\phi_1(x)=2J_1(2t)+xJ_2(2t)+x \cos(2t\sin(\theta))-\frac{x}{2\pi}\int_0^\pi\frac{\sin(2t\sin\varphi)\sin(\varphi)}{\cos(\theta)
-\cos(\varphi)}d\varphi
\end{equation*}
\end{proposition}

\textbf{Proof}.
These equalities follow from Proposition \ref{exp-evol-P} by replacing $k$ by 0 and 1 and using
the summations in Lemma  \ref{Bessel-sum}.
$\qquad\square$\\
\begin{proposition}\label{PQ-mom-02}{\rm
Assuming that $\omega_{n}=1$ for any $n\ge 1$, on the set $\Gamma_0$ one has, for any $t\in\mathbb{C}$,
\begin{equation}\label{PQ-mom-01d1-rep1}
e^{itX}
=\sum_{m_+,m_-,p\ge 0}\frac{(it)^{m_++m_-+2p}}{(m_++m_-+2p)!}
|\Theta_{m_++m_-+2p}(m_+,m_-)|(a^+)^{m_+}(a^-)^{m_-}
\end{equation}
\begin{equation}\label{PQ-mom-01d1-rep2}
	=\sum_{m,n\ge0}i^{m+n}(J_{m+n+2}(2t)+J_{m+n}(2t))(a^+)^m(a^-)^n
\end{equation}
}\end{proposition}
\textbf{Proof}.
The exponential series for $e^{itX}$ is
\begin{align*}
e^{itX}  &  =\sum_{m\ge 0}\frac{( it)^{m}}{m!}X^{m}=\sum_{m\ge 0}\sum_{\varepsilon\in\{-,+\}^{m}}\frac{(it)^{m}}{m!}
a^{\varepsilon(1)}\cdots a^{\varepsilon( m)  }\label{PQ-mom-01a1}
\end{align*}
and one can see that it is obtained from \eqref{itP-2}
replacing $t$ by $it$ and suppressing the factor $(-1)^{\nu_+(\varepsilon)}$, that later in the proof
is identified with $(-1)^{p+m_+}$. With these replacements in the proof of
Proposition \ref{PQ-mom-01}, one arrives to \eqref{PQ-mom-01d1-rep1}.
Finally, using Lemma \ref{df-Imn(t)-01}, one obtains the second identity \eqref{PQ-mom-01d1-rep2}.
$\qquad\square$

\begin{proposition}\label{exp-evol-Q}
	For any $t\in\mathbb{C}$, one has
		\begin{equation*}
		e^{itX}\Phi_0(x)=e^{itx}-ixJ_1(2t)
	\end{equation*}
	and for any $k\ge 1$,
	\begin{equation*}
		e^{itX}\Phi_k(x)=J_{0}(2t)\Phi_{k}(x)+x\sum_{n=1}^{k}i^nJ_{n}(2t)\phi_{k-n+1}(x)+\sum_{n=0}^{k}\sum_{m\ge n+2}i^{m}J_{m}(2t)T_{m+k-2n}(x)
	\end{equation*}
\end{proposition}

\textbf{Proof}.
Replacing the factor $(-1)^{m}$ by $i^{m+n}$ in the proof of Proposition \ref{exp-evol-P}, one arrives to
	\begin{align*}
		e^{itX}\Phi_k(x)=&\sum_{n=0}^{k}\sum_{m\ge0}i^{m+n}(J_{m+n+2}(2t)+J_{m+n}(2t))\Phi_{m+k-n}(x)
		\\=&\sum_{n=0}^{k}\sum_{m\ge n}i^{m}(J_{m+2}(2t)+J_{m}(2t))\Phi_{m+k-2n}(x)
		\\=&\sum_{n=0}^{k}i^nJ_{n}(2t)\Phi_{k-n}(x)+i^{n+1}J_{n+1}(2t)\Phi_{k+1-n}(x)
		\\&+\sum_{m\ge n+2}i^{m}J_{m}(2t)\left\{\Phi_{m+k-2n}(x)-\Phi_{m+k-2n-2}(x) \right\}
		\\=&J_{0}(2t)\Phi_{k}(x)+\sum_{n=1}^{k}i^nJ_{n}(2t)[\Phi_{k-n}(x)+\Phi_{k+2-n}(x)]
		\\&+\sum_{n=0}^{k}\sum_{m\ge n+2}(-1)^{m-n}J_{m}(2t)T_{m+k-2n}(x)
		\\=&J_{0}(2t)\Phi_{k}(x)+x\sum_{n=1}^{k}i^nJ_{n}(2t)\phi_{k-n+1}(x)+\sum_{n=0}^{k}\sum_{m\ge n+2}i^{m}J_{m}(2t)T_{m+k-2n}(x).
	\end{align*}
Now, replacing $k$ by 0 and using \eqref{Jacobi-Anger}, one obtains
\begin{align*}
	e^{itX}\Phi_0(x)&=J_{0}(2t)+\sum_{m\ge 2}i^{m}J_{m}(2t)T_{m}(x)
	\\&=J_{0}(2t)+2\sum_{m\ge 2}i^{m}J_{m}(2t)\cos[m\, arcos(x/2)]
	\\&=e^{itx}-ixJ_1(2t)
\end{align*}
which proves the proposition. $\qquad\square$\\
\begin{proposition}\label{vac-exp-eitP-eitX}{\rm
For any $t\in\mathbb{C}$, the vacuum expectations of $e^{itP}$ and $e^{itX}$ are equal and given by
\begin{equation}\label{vac-exp-eitP-eitQ}
			\langle\Phi_0,e^{itP}\Phi_0\rangle
			=\langle\Phi_0,e^{itX}\Phi_0\rangle
			=\sum_{p=0}^\infty\frac{(-t^2)^{p}}{(2p)!} C_p=\frac{J_1(2t)}{t}
\end{equation}
}\end{proposition}
\textbf{Proof}.
	Putting $m=n=0$ in \eqref{df-Imn(t)}, since (c.f. (3.14) on \cite{[AcHamLu21]})
	$$
	|\Theta_{2p}(0,0)|=C_p
	$$
	one finds,
	$$
	\langle \Phi_0,e^{itP}\Phi_0\rangle
	=\sum_{p\ge 0}\frac{t^{2p}}{(2p)!}(-1)^{p}|\Theta_{2p}(0,0)|
	=\sum_{p\ge 0}\frac{(-t^{2})^{p}}{(2p)!}C_p=\frac{J_1(2t)}{t}
	$$
	which is the second identity in \eqref{vac-exp-eitP-eitQ}. The first one follows from
	from Lemma \ref{X-P-same-vac-distr} because for the semi--circle law the moment problem
	is determined.
$\qquad\square$\\
\begin{remark}
Since the vacuum distributions of $X$ and $P$ coincide with the semi-circle distribution $\mu$,
the vacuum expectations of $e^{itP}$ and $e^{itX}$ coincide with the \textbf{the characteristic function} of the measure $\mu$ which can be computed directly by
expending $e^{ixt}$ into a power series then  using the fact that odd-moments are zero and even-moments are the Catalan numbers interchanging integral and sum, one obtains:
\begin{equation*}
		\varphi_\mu(t)=\sum_{n\ge 0}\frac{(it)^n}{n!}\int_\mathbb{R} x^nd\mu(x)=\sum_{n\ge 0}\frac{(-1)^n}{n+1}\binom{2n}{n}\frac{t^{2n}}{(2n)!}=\frac{J_1(2t)}{t}
	\end{equation*}

\end{remark}

\subsection{The $1$--parameter $*$--automorphism groups associated to $P$}

In this section, we determine the action of $e^{itP}( \ \cdot \ )e^{-itP}$ on the algebra
generated by creation and annihilation operators. For this it is sufficient to determine
this action on $a^+$.
\begin{proposition}{\rm
One has
\begin{equation}\label{P-evol-d-om-Lam}
\partial \omega_{\Lambda,t} := e^{itP} \partial \omega_\Lambda e^{-itP}
=(e^{itP}\Phi_{0})(e^{itP}\Phi_{0})^*
\end{equation}
\begin{equation}\label{P-evol-a+}
a^{+}_t :=e^{itP}a^{+}e^{-itP}
=a^{+} + \omega\int_{0}^{t}ds (e^{isP}\Phi_{0})(e^{isP}\Phi_{0})^*
\end{equation}
\begin{equation}\label{P-evol-X}
X_t :=e^{itP}Xe^{-itP}
= X + 2\omega\int_{0}^{t}ds (e^{isP}\Phi_{0})(e^{isP}\Phi_{0})^*
\end{equation}
}\end{proposition}
\textbf{Proof}.
Since, by definition $\omega_{0}=0$, if $\omega_{\Lambda}$ is given by \eqref{om01},
i.e. $\omega_{\Lambda}=\omega\cdot\mathbf{1}$, one has
\begin{equation}\label{d-omega=0}
\partial\omega_{\Lambda}
=\omega_{\Lambda+1} - \omega_{\Lambda}
= \omega\delta_{0, \Lambda}
= \omega\Phi_{0}\Phi_{0}^*
\end{equation}
From \eqref{df-mom-op} one deduces that $[P,a^{+}]=-i\partial\omega_{\Lambda}$ and this implies
$$
\partial_ta^{+}_t
=ie^{itP}[P,a^{+}]e^{-itP}
= e^{itP}i(-i\partial\omega_{\Lambda})e^{-itP}
=\partial\omega_{\Lambda,t}
$$
or equivalently
$$
a^{+}_t = a^{+} + \int_{0}^{t}ds \partial\omega_{\Lambda,s}
= a^{+} + \int_{0}^{t}ds e^{isP}\partial\omega_{\Lambda}e^{-isP}
$$
$$
\overset{\eqref{d-omega=0}}{=}
a^{+} + \omega\int_{0}^{t}ds e^{isP}\Phi_{0}\Phi_{0}^*e^{-isP}
=a^{+} + \omega\int_{0}^{t}ds (e^{isP}\Phi_{0})(e^{isP}\Phi_{0})^*
$$
Taking adjoint one finds
$$
a_t=a+ \omega\int_{0}^{t}ds e^{isP}\Phi_{0}\Phi_{0}^*e^{-isP}
=a + \omega\int_{0}^{t}ds (e^{isP}\Phi_{0})(e^{isP}\Phi_{0})^*
$$
and summing the two one obtains \eqref{P-evol-X}.
$\qquad\square$\\

\noindent\textbf{Remark}.
From formula \eqref{P-evol-X} one deduced an elegant extension of the action of the usual
quantum mechanical formula $e^{itP}Xe^{-itP}= X+t$.
The main difference is that translation by a number is here replaced by translation by an operator.\\

\noindent In order to give a more explicit form to \eqref{P-evol-a+}, one needs the following result.
\begin{lemma}
		The following inequalities hold
		\begin{equation}\label{est-ser1}
		|I_{0,n}(s)|\le \frac{|s|^{n}}{n!}e^{|s|^2}\left(1+\frac{|s|^2}{2}\right)
		\end{equation}
	and
		\begin{equation}\label{est-eitP}
		|I_{m,n}(s)|\le \frac{|s|^{m+n}}{m!n!}e^{|s|^2}\left(1+\frac{|s|^2}{2}\right)
	\end{equation}
\end{lemma}
\textbf{Proof}.
	One has,
		\begin{align*}
		|J_{n}(2s)|&\le |s|^n\sum_{p\ge0}\frac{|s|^{2p}}{p!(n+p)!}\le \frac{|s|^n}{n!}e^{|s|^2}
	\end{align*}
Then,
	\begin{align*}
		|I_{0,n}(s)|&\le\left|J_{n}(2s)\right|+\left|J_{n+2}(2s)\right|
		\\&\le \frac{|s|^{n}}{n!}e^{|s|^2}\left(1+\frac{|s|^2}{(n+1)(n+2)}\right)
		\\&\le \frac{|s|^{n}}{n!}e^{|s|^2}\left(1+\frac{|s|^2}{2}\right)
	\end{align*}
Now, to prove \eqref{est-eitP}, notice that
\begin{equation*}
	(n+m)! = \frac{(n+m)!}{n!m!}n!m! \ge n!m!
\end{equation*}
Therefore,
$$
		|I_{m,n}(s)|=	|I_{0,m+n}(s)|\le\frac{|s|^{m+n}}{(m+n)!}e^{|s|^2}\left(1+\frac{|s|^2}{2}\right)
\le \frac{|s|^{m+n}}{m!n!}e^{|s|^2}\left(1+\frac{|s|^2}{2}\right)
$$
which proves the Lemma. $\qquad\square$\\
\begin{theorem}\label{th:expl-espr-a+t}{\rm
One has, for each $t\in \mathbb{R}$,
\begin{equation}\label{expl-espr-a+t}
a^{+}_t = a^{+}  +\omega\sum_{m,n,p,q\ge 0}\frac{(-1)^{m+n+p+q}(m+1)(n+1)}{p!q!(p+m+1)!(q+n+1)!} \frac{t^{m+n+2p+2q+1}}{m+n+2p+2q+1}\Phi_m\Phi_{n}^*
\end{equation}
}\end{theorem}
\textbf{Proof}.
Using \eqref{P-evol-a+}, one has
\begin{align*}
	a^{+}_t =e^{itP}a^{+}e^{-itP} \overset{\eqref{eitP-phi0}}{=}a^{+}  + \omega\int_{0}^{t}ds
	\left(\sum_{m\ge 0}I_{m,0}(s)\Phi_m\right)\left(\sum_{n\ge 0}I_{n,0}(s)\Phi_{n}\right)^*
\end{align*}

$$
=
a^{+}  + \omega\int_{0}^{t}ds
\left(\sum_{m\ge 0}\sum_{p\ge 0}\frac{s^{m+2p}}{(m+2p)!}(-1)^{p+m}|\Theta_{m+2p}(m,0)|\Phi_m
\right)
$$
\begin{equation}\label{int-ser1}
\left(\sum_{n\ge 0}\sum_{q\ge 0}\frac{s^{n+2q}}{(n+2q)!}(-1)^{q+n}|\Theta_{n+2q}(n,0)|\Phi_{n}\right)^*
\end{equation}
where \eqref{est-eitP} guarantees that the series in \eqref{int-ser1} is absolutely convergent.
Therefore one can exchange series and integral obtaining
$$
a^{+}  + \sum_{m,n,p,q\ge 0}
\frac{|\Theta_{m+2p}(m,0)|(-1)^{q+n+p+m}|\Theta_{n+2q}(n,0)|}{(m+2p)!(n+2q)!}
\Phi_m\Phi_{n}^* \omega\int_{0}^{t}ds s^{m+2p+n+2q}
$$
\begin{multline}
	= a^{+}
	+ \omega\sum_{m,n,p,q\ge 0}
	\frac{|\Theta_{m+2p}(m,0)|(-1)^{q+n+p+m}|\Theta_{n+2q}(n,0)|}{(m+2p)!(n+2q)!}\\
	\frac{t^{m+2p+n+2q+1}}{m+2p+n+2q+1}\Phi_m\Phi_{n}^*
\end{multline}

Finally, using the expression (3.17) of \cite{[AcHamLu21]}, one finds \eqref{expl-espr-a+t}.
$\qquad\square$\\

\begin{remark}
	Using the expression \eqref{I0n(t)}, one can also rewrite \eqref{expl-espr-a+t} as
	\begin{equation}\label{expl-espr-a+t1}
		a^{+}_t = a^{+}  +\omega\sum_{m,n\ge 0}(-1)^{m+n}(m+1)(n+1)\int_{0}^{t}\frac{ds}{s^2}J_{m+1}(2s)J_{n+1}(2s)\Phi_m\Phi_{n}^*
	\end{equation}
\end{remark}

\section{Free evolutions associated to the kinetic energy operator}\label{sec:evol-eitP2}

\subsection{The free Hamiltonian group $e^{itP^2}$}

\begin{proposition}\label{p-sqrt01}{\rm
If $\omega_{n}=1$ for any $n\geq1$,
for any $t$,
\begin{equation}\label{p-sqrt01a}
e^{itP^{2}}
=\sum_{m,n\geq0}I_{m,n}^{\left(  2\right)  }(t)(a^{+})^{m}a^{n}
\end{equation}
where 
for any $m,n\in\mathbb{N}$,
\begin{equation}\label{p-sqrt01b}
I_{m,n}^{\left(  2\right)  }(t)
=\chi_{2\mathbb{N}}\left(  m+n\right)\left(-1\right)^{\frac{3m+n}{2}}  \sum_{p\geq0}
\frac{(it)^{\frac{m+n}{2}+p}}{\left(  \frac{m+n}{2}+p\right)  !}
|\Theta_{m+n+2p}(m,n)|
\end{equation}
}\end{proposition}
\textbf{Proof}.
Expanding $e^{itP^{2}}$ we obtain
$$
e^{itP^{2}}=\sum_{m\geq0}\frac{\left(  it\right)  ^{m}}{m!}P^{2m}
\overset{P=-i\left(  a-a^{+}\right)  }{=}
\sum_{m\geq0}\frac{\left(-it\right)  ^{m}}{m!}\sum_{\varepsilon\in\{-,+\}^{2m}}(-1)^{\nu_{+}
(\varepsilon)}a^{\varepsilon(1)}\cdots a^{\varepsilon(2m)}
$$
where $\varepsilon\in\{-,+\}^{2m}$ determines uniquely $m_{\pm}$
such that $2m-m_{+}-m_{-}$ is even and is denoted by $2p$. Since
$\omega_{n}\equiv1$ for any $n\geq1$,
\begin{align*}
a^{\varepsilon(1)}\cdots a^{\varepsilon(2m)}  & =\left(  a^{+}\right)
^{m_{+}}a^{m_{-}}\\
\nu_{+}(\varepsilon)  & =m_{+}+\frac{2m-m_{+}-m_{-}}{2}=m+\frac{m_{+}-m_{-}%
}{2}=m_{+}+p
\end{align*}
Since $2m-m_{+}-m_{-}$ is even if and only if $m_{+}+m_{-}$ is even,
formula (3.8) of \cite{[AcHamLu21]}
with $c=1$, for the solution of the inverse normal order problem, gives
\begin{align*}
& \sum_{m\geq0}\frac{\left(-it\right)  ^{m}}{m!}\sum_{\varepsilon
\in\{-,+\}^{2m}}(-1)^{\nu_{+}(\varepsilon)}a^{\varepsilon(1)}\cdots a^{\varepsilon(2m)}\\
& =\sum_{m_{\pm},p\geq0,\ m_{+}+m_{-}\text{ even}}\frac{\left(-it\right)
^{\frac{m_{+}+m_{-}}{2}+p}}{\left(\frac{m_{+}+m_{-}}{2}+p\right)  !}\left(
-1\right)  ^{m_{+}+p}|\Theta_{m_{+}+m_{-}+2p}(m_{+},m_{-})|\left(
a^{+}\right)  ^{m_{+}}a  ^{m_{-}}\\
& =\sum_{m,n\geq0,\ m+n\text{ even}}(-1)  ^{\frac{3m+n}{2}} \sum_{p\geq0}\frac{(it)^{\frac{m+n}{2}+p}%
}{\left(\frac{m+n}{2}+p\right)  !}%
|\Theta_{m+n+2p}(m,n)|\left(a^{+}\right)  ^{m}a  ^{n}\\
& =\sum_{m,n\geq0}\chi_{2\mathbb{N}}\left(m+n\right) (-1)
^{\frac{3m+n}{2}} \sum_{p\geq0}%
\frac{(it)^{\frac{m+n}{2}+p}}{\left(\frac{m+n}{2}+p\right)  !}|\Theta_{m+n+2p}(m,n)|\left(a^{+}\right)  ^{m}a  ^{n}
\end{align*}
$\hfill\qquad\square$

\begin{proposition}
	For any $m,n\in\mathbb{N}$,
		\begin{equation*}
		I_{m,n}^{\left(2\right)  }(t)	
=(-1)^mI_{0,m+n}^{\left(2\right)  }(t)=(-1)^nI_{m+n,0}^{\left(2\right)}(t)
	\end{equation*}
with
	\begin{equation}\label{p-sqrt01c}
		I_{0,n}^{\left(  2\right)  }(t)
	=\chi_{2\mathbb{N}}\left(  n\right)
\frac{(-it)^{\frac{n}{2}}}{\left(  \frac{n}{2}\right)!}\,  {}_1F_1\left(\frac{n+1}{2}; n+2;4it\right)
	\end{equation}
where ${}_1F_1$ is the confluent hypergeometric function.
\end{proposition}
\textbf{Proof}.
Replacing $|\Theta_{m+n+2p}(m,n)|$ in  \eqref{p-sqrt01b} by its expression given by formula
(3.17) of \cite{[AcHamLu21]},
one finds
	\begin{align*}
			I_{m,n}^{\left(  2\right)  }(t)	&= \chi_{2\mathbb{N}}\left(  m+n\right)\left(-1\right)^{\frac{3m+n}{2}}  \sum_{p\geq0}
			\frac{(it)^{\frac{m+n}{2}+p}}{\left(  \frac{m+n}{2}+p\right)  !}  \frac{m+n+1}{2p+m+n+1}\frac{(2p+m+n+1)!}{p!(p+m+n+1)!}
			\\&= \chi_{2\mathbb{N}}\left(  m+n\right)\left(-1\right)^{\frac{3m+n}{2}} (m+n+1) \sum_{p\geq0}
			\frac{(it)^{\frac{m+n}{2}+p}}{\left(  \frac{m+n}{2}+p\right)  !}  \frac{(2p+m+n)!}{p!(p+m+n+1)!}
	\end{align*}
Next, using the Legendre duplication formula (see e.g.  \cite[p. 24]{Rainville}), one obtains
	\begin{align*}
	I_{m,n}^{\left(  2\right)  }(t)	&= \chi_{2\mathbb{N}}\left(  m+n\right)\left(-1\right)^{\frac{3m+n}{2}} (m+n+1) \sum_{p\geq0}
	\frac{(4it)^{\frac{m+n}{2}+p}}{\sqrt{\pi} p!}  \frac{\Gamma(p+\frac{m+n+1}{2})}{\Gamma(p+m+n+2)}
\end{align*}
or equivalently
\begin{multline*}
		I_{m,n}^{\left(  2\right)  }(t)
= \chi_{2\mathbb{N}}\left(  m+n\right)\left(-1\right)^{m} (m+n+1)
\frac{	(-4it)^{\frac{m+n}{2}}\Gamma(\frac{m+n+1}{2})}{\sqrt{\pi}\Gamma(m+n+2)}\\
\times{}_1F_1\left(\frac{m+n+1}{2}; m+n+2; 4it\right)
\end{multline*}
Using again the Legendre duplication formula, one gets
\begin{align*}
	I_{m,n}^{\left(  2\right)  }(t)	
&= \chi_{2\mathbb{N}}\left(  m+n\right)\left(-1\right)^{m}
\frac{	(-it)^{\frac{m+n}{2}}}{\Gamma(1+\frac{m+n}{2})}
\,{}_1F_1\left(\frac{m+n+1}{2}; m+n+2; 4it\right)
\end{align*}

\begin{remark}
	Since
	\begin{equation*}
		\chi_{2\mathbb{N}}\left(  m+n\right) =\chi_{2\mathbb{N}}
\left(  m\right) \chi_{2\mathbb{N}}\left(  n\right) +\chi_{2\mathbb{N}+1}\left(  m\right) \chi_{2\mathbb{N}+1}\left(  n\right)
	\end{equation*}
	the equality \eqref{p-sqrt01c} can also be expressed as
	\begin{equation*}
			I_{m,n}^{\left(  2\right)  }(t)	=
		\begin{cases}
				\frac{(-it)^{j}}{ j!}\,  {}_1F_1\left(j+\frac{1}{2}; 2(j+\frac{1}{2})+1; 4it\right),& m=2k, n=2l\\
				-	\frac{(-it)^{j+1}}{\left(  j+1\right)!}\,  {}_1F_1\left(j+\frac{3}{2}; 2(j+\frac{3}{2})+1; 4it\right),& m=2k+1, n=2l+1\\
					0,& otherwise
		\end{cases}
	\end{equation*}
where $j=k+l$.
\end{remark}

\subsection{The $1$--parameter $*$--automorphism groups associated to $P^2$}

With the notation $\partial \omega_\Lambda:=\omega_{\Lambda+1}-\omega_\Lambda$, denote
\begin{equation}\label{P2-evol-dom-a+}
\partial \omega_{\Lambda,t} := e^{itP^{2}} \partial \omega_\Lambda e^{-itP^{2}}
=(e^{itP^{2}}\Phi_{0})(e^{itP^{2}}\Phi_{0})^*
\ ; \
a^{+}_t :=e^{itP^{2}}a^{+}e^{-itP^{2}}
\end{equation}
From Theorem 1 in \cite{[AcEllLu20-mom]} we know that, for general $1$MIFS, the integral forms
of the equation of motion for $a^{+}$ and $\partial\omega_{\Lambda, t}$ are respectively
\begin{equation}\label{ev-a+t-p2-int1}
a^{+}_t = a^{+} +   P \circ \int_{0}^{t}ds \ (\partial\omega_{\Lambda})_{s}
\end{equation}
\begin{multline}\label{ev-domegat-p2-int1}
(\partial\omega_{\Lambda})_{t} = \partial\omega_{\Lambda} + P \circ\left(a^{+}\circ \int_{0}^{t}ds (\partial^2\omega_{\Lambda})_{s}\right)\\
+ P \circ\left(P \circ \int_{0}^{t}ds  \int_{0}^{s}dr\ (\partial\omega_{\Lambda})_{r}\circ(\partial^2\omega_{\Lambda})_{s}\right)
\end{multline}

where for two linear operators $A$ and $B$, $A\circ B := AB + B^*A^*$
(see \cite{[AcEllLu20-mom]}).
\begin{remark}
Notice that
$$
\partial_t a^{+}_t =e^{itP^{2}}[P^2,a^{+}]e^{-itP^{2}}
=e^{itP^{2}}(P\circ \Phi_0\Phi_0^*)e^{-itP^{2}}
=\left(ie^{itP^{2}}\Phi_1\right)\circ\left(e^{-itP^{2}}\Phi_0\right)^*
$$
Hence, one can rewrite \eqref{ev-a+t-p2-int1} as
\begin{equation}
		a^{+}_t = a^{+} +   \int_{0}^{t}ds \ \left(ie^{isP^{2}}\Phi_1\right)\circ\left(e^{-isP^{2}}\Phi_0\right)^*
	\end{equation}
\end{remark}

\begin{theorem}\label{ev-a+t-P2-expl}{\rm
One has
\begin{equation}\label{expl-evol-a+}
a^{+}_t = a^{+} +    \sum_{m, p, n, q\geq0}
\frac{(-1)^{m+q}|\Theta_{2m+2p}(2m,0)||\Theta_{2n+2q}(2n,0)|}
{\left(m+p\right)!\left(n+q\right)!}
\end{equation}
$$
\frac{(it)^{m+n+p+q+ 1}}{m+n+p+q+ 1}\left(T_{2m-1}\Phi_{2n}^*- \Phi_{2n}T_{2m-1}^* \right)
$$
where we recall from the formula (3.17) of \cite{[AcHamLu21]} that
\begin{equation*}
	|\Theta_{2m+2p}(2n,0)|
=\frac{2m+1}{2p+2m+1}\binom{2p+2m+1}{p}
=\frac{(2m+1)(2m+2p)!}{p!(2m+p+1)!}.
\end{equation*}
}\end{theorem}
\textbf{Proof}.
One has
\begin{equation}\label{ev-a+t-p2-int2}
P \circ \int_{0}^{t}ds \ (\partial\omega_{\Lambda})_{s}
\overset{\eqref{P2-evol-dom-a+}}{=}
P \circ \int_{0}^{t}ds \ (e^{isP^{2}}\Phi_{0})(e^{isP^{2}}\Phi_{0})^*
\end{equation}
\eqref{p-sqrt01a} and \eqref{p-sqrt01b} imply that
$$
e^{itP^{2}}\Phi_0
=\sum_{m,p\geq0}\chi_{2\mathbb{N}}(m)
\frac{(-1)^{\frac{3m}{2}}(it)^{\frac{m}{2}+p}}{\left(\frac{m}{2}+p\right)!}
|\Theta_{m+2p}(m,0)|\Phi_m
$$
\begin{equation}\label{e(itP2)Phi0}
=\sum_{m,p\geq0}(-1)^m
\frac{(it)^{m+p}|\Theta_{2m+2p}(2m,0)|}{(m+p)!} \Phi_{2m}
\end{equation}
Replacing $e^{isP^{2}}\Phi_{0}$ in \eqref{ev-a+t-p2-int2} by its expression given by
\eqref{e(itP2)Phi0}, one finds, using the fact that the coefficients in the series in
\eqref{e(itP2)Phi0} are real,
$$
P \circ \int_{0}^{t}ds \ (\partial\omega_{\Lambda})_{s}
\overset{\eqref{P2-evol-dom-a+}}{=}
P \circ \int_{0}^{t}ds \ (e^{isP^{2}}\Phi_{0})
(e^{isP^{2}}\Phi_{0})^*
$$
$$
\overset{\eqref{e(itP2)Phi0}}{=}
P \circ \int_{0}^{t}ds \
\sum_{m, p\geq0}(-1)^m
\frac{(is)^{m+p}|\Theta_{2m+2p}(2m,0)|}{(m+p)!} \Phi_{2m}
$$
$$
\sum_{n, q\geq0}(-1)^n
\frac{(-is)^{n+q}|\Theta_{2n+2q}(2n,0)|}{(n+q)!} \Phi_{2n}^*
$$
$$
= \sum_{m, p, n, q\geq0}
\frac{(-1)^{m+q}|\Theta_{2m+2p}(2m,0)||\Theta_{2n+2q}(2n,0)|}
{\left(m+p\right)!\left(n+q\right)!}
$$
$$
\int_{0}^{t}ds \ s^{m+n+p+q} P \circ \Phi_{2m}\Phi_{2n}^*
$$
$$
= \sum_{m, p, n, q\geq0}
\frac{(-1)^{m+q}|\Theta_{2m+2p}(2m,0)||\Theta_{2n+2q}(2n,0)|}
{i\left(m+p\right)!\left(n+q\right)!}
\frac{(it)^{m+n+p+q+ 1}}{m+n+p+q+ 1}
$$
$$
\left((i(a^{+}-a) \Phi_{2m})\Phi_{2n}^* + \Phi_{2n}(i(a^{+}-a)\Phi_{2m})^* \right)
$$
$$
= \sum_{m, p, n, q\geq0}
\frac{(-1)^{m+q}|\Theta_{2m+2p}(2m,0)||\Theta_{2n+2q}(2n,0)|}
{\left(m+p\right)!\left(n+q\right)!}
\frac{(it)^{m+n+p+q+ 1}}{m+n+p+q+ 1}
$$
$$
\left((\Phi_{2m+1}-\Phi_{2m-1})\Phi_{2n}^*
- \Phi_{2n}(\Phi_{2m+1}-\Phi_{2m-1})^* \right)
$$
$$
\overset{\eqref{connection2}}{=} \sum_{m, p, n, q\geq0}
\frac{(-1)^{m+q}|\Theta_{2m+2p}(2m,0)||\Theta_{2n+2q}(2n,0)|}
{\left(m+p\right)!\left(n+q\right)!}
\frac{(it)^{m+n+p+q+ 1}}{m+n+p+q+ 1}
$$
$$
\left(T_{2m-1}\Phi_{2n}^*
- \Phi_{2n}T_{2m-1}^* \right)
$$


\begin{thebibliography}{99}





\bibitem{[AcBaLuRha15]}
Luigi Accardi, Abdessatar Barhoumi, Yun Gang Lu, Mohamed Rhaima:\\
$*$--Lie algebras canonically associated to Probability Measures on $\mathbb{R}$ with all moments,\\
in: Proceedings of the XI International Workshop ''Lie Theory and Its Applications
in Physics'', (Varna, Bulgaria, June 2015), Springer Proceedings in Mathematics and
Statistics Vol. 191, ed. V. Dobrev, Springer (2016) 3--21\\


\bibitem{[AcBo98]}
Accardi L., Bozejko M.:\\
Interacting Fock spaces and Gaussianization of probability measures,\\
Infin. Dimens. Anal. Quantum Probab. Relat. Top. (IDA-QP) 1 (4) (1998) 663-670\\
Volterra Preprint N. 321 (1998)


\bibitem{[AcEllLu20-mom]}
Luigi Accardi, Abdon Ebang Ella, Yun Gang Lu:\\
Quantum symmetries and momentum operators associated to classical random variables,\\
In preparation (2020)

\bibitem{[AcHamLu21]}
Luigi Accardi, Tarek Hamdi, Yun Gang Lu:\\
The quantum mechanics canonically associated to free probability.
Part II: The normal and inverse normal order problem.\\
submitted for publication 2021


\bibitem{[AcLuVo97a-QED-Hilb-mod]}
Accardi L., Lu Y.G., Volovich I.:\\
The QED Hilbert module and Interacting Fock spaces,\\
Publications of IIAS (Kyoto), N1997-008 (1997)

\bibitem{[APS96]}
K. Astala, L. Paivarinta and E. Saksman:\\
The finite Hilbert transform in weighted spaces, \\
Proceedings of the Royal Society of Edinburgh: Mathematics 126 (1996) 1157--1167

\bibitem{[BuNe71]}
P. L. Butzer and R. J. Nessel:\\
Fourier Analysis and Approximation,\\
Vol. 1 Birkh\"auser-Verlag, Basel (1971)

\bibitem{CPVWJ08}
 Annie Cuyt, Vigdis Petersen, Brigitte Verdonk, Haakon Waadeland et William B. Jones:\\
  Handbook of Continued Fractions for Special Functions,\\
   Springer, (2008)

\bibitem{[HoOb06]}
A. Hora and N. Obata,
Quantum Probability and Spectral Analysis of Graphs,
Theoretical and Mathematical Physics, Springer (2007)

\bibitem{[Kap1905]}
W. Kapteyn:\\
The values of some definite integrale connected with Bessel functions,\\
 KNAW, Proceedings, 7, 1904-1905, Amsterdam, (1905) 375--376

\bibitem{[KaTo12]}
A. Katsevich, A. Tovbis:\\
Finite Hilbert transform with incomplete data: null-space and singular values,\\
Inverse  Probl. 28 (2012) 105006

\bibitem{[King09]}
F. King:\\
Hilbert Transforms,\\
Encyclopedia of Mathematics and its applications,
Vol.1, 124, Cambridge (2009)


\bibitem{[MaHa03]}
J. C. Mason and D. C. Handscomb: \\
Chebyshev Polynomials,\\
CRC Press (2003)

\bibitem{Rainville}
E. D. Rainville:\\
Special functions,\\
The Macmillan Co. New York (1960)

\bibitem{[Tricomi57]}
F. G. Tricomi:\\
Integral Equations,\\
Interscience Publisher Inc (1957)


\end{thebibliography}
\end{document}